\documentclass[letterpaper,10pt]{amsart}
\usepackage{amssymb}
\usepackage{amsthm}
\usepackage[utf8]{inputenc}
\usepackage[margin=1.0in]{geometry}
\usepackage{pst-node}
\usepackage{tikz-cd} 
\usepackage{thmtools}
\usepackage{thm-restate}
\usepackage{hyperref}
\usepackage{amsrefs}
\usepackage[noabbrev,capitalize]{cleveref}
\usepackage{amsmath}
\usepackage{amsfonts}
\usepackage{tikz-cd}
\usepackage{mathtools}
\usepackage{bbm}
\usepackage{hyperref}
\usepackage{blindtext, xcolor}
\usepackage{comment}
\usepackage{romannum}
\usepackage{multirow}
\usepackage{array}
\usepackage{wrapfig}
\usepackage{graphicx}
\usepackage{caption}

\DeclareMathOperator{\Hom}{Hom}

\DeclareMathOperator{\Stab}{Stab}

\DeclareMathOperator{\Spec}{Spec}
\DeclareMathOperator{\Pic}{Pic}

\DeclareMathOperator{\Len}{Len}

\DeclareMathOperator{\perc}{\%}

\DeclareMathOperator{\lv}{\lvert}
\DeclareMathOperator{\rv}{\rvert}
\DeclareMathOperator{\Spe}{Spectrum}
\DeclareMathOperator{\Disc}{Disc}

\DeclareMathOperator{\Cone}{Cone}

\newcommand{\ZZ}{\mathbb{Z}}      
\newcommand{\QQ}{\mathbb{Q}}      
\newcommand{\RR}{\mathbb{R}}      
\newcommand{\AAA}{\mathbb{A}}      
\newcommand{\cO}{\mathcal{O}}      
\newcommand{\CC}{\mathbb{C}}      
\newcommand{\PP}{\mathbb{P}}      

\newcommand{\mfa}{\mathfrak{a}}

\newcommand{\mfT}{\mathfrak{T}}

\newcommand{\la}{\langle}    
\newcommand{\ra}{\rangle}

\newcommand{\cF}{\mathcal{F}}

\newcommand{\cM}{\mathcal{M}}    
\newcommand{\cL}{\mathcal{L}}

\newcommand{\Mod}[1]{\ (\mathrm{mod}\ #1)}

\newtheorem{theorem}{Theorem}[section]

\newtheorem{lemma}[theorem]{Lemma}
\newtheorem{proposition}[theorem]{Proposition}
\newtheorem{corollary}[theorem]{Corollary}

\theoremstyle{definition}
\newtheorem{conjecture}[theorem]{Conjecture}
\newtheorem{example}[theorem]{Example}

\newtheorem{question}[theorem]{Question}
\newtheorem{observation}[theorem]{Observation}

\theoremstyle{remark}

\newtheorem{definition}[theorem]{Definition}

\usepackage[english]{babel}
\usepackage[utf8]{inputenc}

\title{Bounds on Successive Minima of Orders in Number Fields and Scrollar Invariants of Curves}

\author{Sameera Vemulapalli}
\address{Department of Mathematics, Harvard University}
\email{vemulapalli@math.harvard.edu}


\usepackage[english]{babel}
\usepackage[utf8]{inputenc}

\begin{document}

\maketitle
\pagenumbering{arabic}

\begin{abstract}
Orders and fractional ideals in number fields provide interesting examples of lattices. We ask: what lattices arise from orders in number fields? We prove that all nontrivial  multiplicative constraints on successive minima of orders come from multiplication. Moreover, inspired by a conjecture of Lenstra, for infinitely many positive integers $n$ (including all $n < 18$), we explicitly determine all multiplicative constraints on successive minima of orders in degree $n$ number fields. We also prove analogous results for scrollar invariants of curves. 
\end{abstract}

\tableofcontents

\section{Introduction}
Orders and ideals in number fields of degree $n$ provide interesting examples of lattices via their natural embeddings into $\RR^n$ using their real and complex places. The shapes of these lattices are constrained due to multiplication: the length of the product of two vectors is roughly bounded above by the product of the lengths. By studying this multiplicative structure, we make these constraints explicit. 

More precisely, let $\mfa$ be a fractional ideal of an order $\cO$ in a degree $n$ number field $K$. Denote the nonzero homomorphisms of $K$ into $\CC$ by $\sigma_1,\dots,\sigma_{n}$, and define 
\[
	\lvert x \rvert \coloneqq \sqrt{\frac{1}{n}\sum_{i = 1}^{n} \lvert\sigma_i(x)\rvert^2}
\]
for $x \in K$. Set $[n] \coloneqq \{0,\dots,n-1\}$. For $i \in [n]$, let $\lambda_i(\mfa)$ be the $i$th successive minima of $\mfa$ with respect to this norm, e.g., the smallest positive real number $r$ such that $\mfa$ contains at least $i+1$ linearly independent elements of length $\leq r$. 

\begin{theorem}[Bhargava, Lenstra, unpublished]
\label{thm:bl}
If $K$ has no nontrivial proper subfields, and $\mfa_1,\mfa_2,\mfa_3$ are fractional ideals such that $\mfa_1\mfa_2 = \mfa_3$, then
\[
    \lambda_{i + j}(\mfa_3) \leq \sqrt{n} \lambda_i(\mfa_1)\lambda_j(\mfa_2)
\]
for any integers $0 \leq i,j \leq i + j < n$.
\end{theorem}

The assumption that $K$ has no nontrivial proper subfields is necessary for \Cref{thm:bl}. Take for example, the order $\cO = \ZZ[i,\sqrt{101}]$ and take $\mfa_1 = \mfa_2 = \mfa_3 = \cO$. Then $\lambda_1(\cO) = \lvert i \rvert = 1$ and $\lambda_2(\cO) = \lvert \sqrt{101} \rvert = \sqrt{101}$, so
\[
    \lambda_2(\cO) > \sqrt{4}\lambda_1(\cO)\lambda_1(\cO).
\]
Allowing for the existence of subfields, we have a generalization of \Cref{thm:bl} (indeed, \Cref{thm:bl} is a corollary of \Cref{thm:bound-subfield}). For positive integers $i,j$, let $i\%j$ denote the remainder when dividing $i$ by $j$.

\begin{theorem}
\label{thm:bound-subfield}
Fix integers $0 \leq i,j \leq i + j < n$. Suppose that for every integer $m$ such that $K$ has a degree $m$ subfield, we have $(i\%m) + (j\%m) = (i + j)\%m$. Then for any three fractional ideals $\mfa_1,\mfa_2,\mfa_3$ with $\mfa_1\mfa_2 = \mfa_3$, we have
\[
    \lambda_{i + j}(\mfa_3) \leq \sqrt{n}\lambda_i(\mfa_1)\lambda_j(\mfa_2).
\]
\end{theorem}

We now say a few words illustrating the key idea behind the proof of \Cref{thm:bound-subfield}. Let $v_0,\dots,v_{n-1}$ be a set of linearly independent vectors in $\mfa_1$ with the property that $\lvert v_i \rvert = \lambda_i(\mfa_1)$. Similarly, let $u_0,\dots,u_{n-1}$ be a set of linearly independent vectors in $\mfa_2$ with the property that $\lvert u_i \rvert = \lambda_i(\mfa_2)$.

Given a field extension $K/L$ and two $L$-vector spaces $I,J \subseteq K$, set $IJ \coloneqq \{vu : v \in I, \; u \in J\}$, where the multiplication is simply multiplication in the field $K$. Given elements $v_1,\dots,v_{\ell} \in K$, let $L\langle v_1,\dots,v_{\ell}\rangle$ denote the $L$-vector space spanned by the $v_i$. The crucial tool in the proof of \Cref{thm:bound-subfield} is the following proposition. 

\begin{proposition}
\label{prop:mult-bound}
Fix integers $0 \leq i,j < n$. Set $k \coloneqq \dim_{\QQ}\QQ\langle v_0,\dots,v_i \rangle \QQ\langle u_0,\dots,u_j\rangle -1$. Then
\[
    \lambda_k(\mfa_3) \leq \sqrt{n}\lambda_i(\mfa_1)\lambda_j(\mfa_2).
\]
\end{proposition}

To use \Cref{prop:mult-bound} to prove \Cref{thm:bound-subfield}, we prove lower bounds on the dimension of the product space $\dim_{\QQ}\QQ\langle x_0,\dots,x_i \rangle \QQ\langle y_0,\dots,y_j\rangle$ using theorems from additive combinatorics. To illustrate this approach in an elementary case, we do an example.

\begin{example}
Let $\cO$ be an order in a cubic field and set $\mfa_1 = \mfa_2 = \mfa_3 = \cO$. As above, let $v_0,v_1,v_2 \in \cO$ be linearly independent elements such that $\lambda_i(\cO) = \lvert v_i \rvert$. Without loss of generality, we may take $v_0 = 1$. Set $i = j = 1$. Then the product space
\[
    \QQ \langle 1,v_1\rangle\QQ \langle 1,v_1\rangle
\]
has dimension $3$; it contains the three linearly independent vectors $\{1,v_1,v_1^2\}$. Therefore, \Cref{prop:mult-bound} implies that
\[
    \lambda_2(\cO) \leq \sqrt{3} \lambda_1(\cO)\lambda_1(\cO). 
\]
\end{example}

\subsection{Bounds on successive minima of orders in number fields}

We now restrict our focus from the successive minima of fractional ideals to the successive minima of orders, e.g., we specialize to the case $\mfa_1 = \mfa_2 = \mfa_3 = \cO$. In this case, $\lambda_0(\cO) = 1$ (see \Cref{lem:lambda-zero}). We ask: as we range across orders $\cO$ in degree $n$ number fields, what are the possible values of the tuples
\[
    (\lambda_1(\cO),\dots,\lambda_{n-1}(\cO)) \in \RR^{n-1}.
\]
It turns out that there are interesting relationships between the successive minima which are not captured by \Cref{thm:bound-subfield}.

\begin{example}
\label{ex:6}
Set $n = 6$. There exist orders in sextic fields with $\lambda_2 > \sqrt{6}\lambda_1\lambda_1$; take for example $\cO = \ZZ[\sqrt{2},\sqrt[3]{101}]$. Similarly, there exist orders in sextic fields with $\lambda_3 > \sqrt{6}\lambda_1\lambda_2$; take $\cO = \ZZ[\sqrt[3]{2},\sqrt{101}]$. However, there do \emph{not} exist orders in sextic fields such that $\lambda_2 > \sqrt{6}\lambda_1\lambda_1$ \emph{and} $\lambda_3 > \sqrt{6}\lambda_1\lambda_2$, as we show below.

Let $\cO$ be an order in a degree $6$ number field and suppose that $\lambda_2 > \sqrt{6}\lambda_1\lambda_1$. Let $1=x_0,\dots,x_5 \in \cO$ be elements such that $\lvert x_i\rvert = \lambda_i(\cO)$. Then \Cref{prop:mult-bound} implies that $\dim_{\QQ}(\QQ \la 1,x_1\ra \QQ \la 1,x_1\ra) \leq 2$, and so $M \coloneqq \QQ \la 1,x_1\ra$ is a quadratic field. So, the product space $\QQ \la 1,x_1\ra \QQ \la 1,x_1,x_2\ra$ is a vector space over the quadratic field $M$. Therefore $\dim_{\QQ}\QQ \la 1,x_1\ra \QQ \la 1,x_1,x_2\ra \geq 4$. Hence, \Cref{prop:mult-bound} implies that $\lambda_3 \leq \sqrt{6}\lambda_1 \lambda_2$.
\end{example}

The contribution of \Cref{thm:containment} is to capture which constraints among successive minima hold jointly. In order to phrase our theorem, we will need the following notation.

\begin{definition}
\label{def:tower-type}
A \emph{tower type} is a $t$-tuple of integers $(n_1,\dots,n_t) \in \ZZ_{>1}^t$ for some $t \geq 1$. We say $\prod_{i = 1}^tn_i$ is the \emph{degree} of the tower type and $t$ is the \emph{length} of the tower type. 
\end{definition}

\noindent Throughout this article, the variable $\mfT$ will refer to a tower type of length $t$ and degree $n$. 

\begin{definition}
\label{def:mixed-radix}
Choose a tower type $\mfT = (n_1,\dots,n_t)$ and $i \in [n]$. Writing $i$ in \emph{mixed radix notation with respect to $\mfT$} means writing
\[
	i = i_1 + i_2 n_1 + i_3 (n_1n_2) + \dots + i_t(n_1\dots n_{t-1})
\]
where $i_s$ is an integer such that $0 \leq i_s < n_s$ for $1\leq s \leq t$. Note that the integers $i_s$ are uniquely determined.
\end{definition}

\begin{definition}
\label{def:overflow}
Fix a tower type $\mfT=(n_1,\dots,n_t)$ and integers $0 \leq i, j \leq i + j < n$. Write $i$, $j$, and $k=i+j$ in mixed radix notation with respect to $\mfT$ as 
\[
	i = i_1 + i_2 n_1 + i_3 (n_1n_2) + \dots + i_t(n_1\dots n_{t-1})
\]
\[
	j = j_1 + j_2 n_1 + i_3 (n_1n_2) + \dots + j_t(n_1\dots n_{t-1})
\]
\[
	k = k_1 + k_2 n_1 + k_3 (n_1n_2) + \dots + k_t(n_1\dots n_{t-1}).
\]
We say the addition $i+j$ \emph{does not overflow modulo $\mfT$} if $i_s + j_s = k_s$ for all $1 \leq s \leq t$. Otherwise, we say the addition $i+j$ \emph{overflows modulo $\mfT$}.
\end{definition}

\begin{theorem}
\label{thm:containment}
Suppose $n$ is a prime power, a product of $2$ primes, or equal to $12$. Let $\cO$ be an order in a degree $n$ number field. Then there exists a tower type $\mfT$, depending only on $\cO$, such that for all $0 \leq i,j \leq i+j < n$, if $i + j$ does not overflow modulo $\mfT$, then
\[
    \lambda_{i + j}(\cO) \ll_n \lambda_i(\cO) \lambda_j(\cO).
\]
\end{theorem}

For every $n$ which is not a prime power, a product of $2$ primes, or equal to $12$, the statement of \Cref{thm:containment} is false; see \Cref{thm:lenstra-spectrum}. Namely, upon fixing such an integer $n$, for every positive real number $c$ there exists an order $\cO$ in a degree $n$ number field such that for every tower type $\mfT$, there exists $0 \leq i,j \leq i+j < n$ such that $i + j$ does not overflow modulo $\mfT$ and 
\[
    \lambda_{i + j}(\cO) > c\lambda_i(\cO) \lambda_j(\cO).
\]

\subsection{The successive minima spectrum}
\Cref{thm:bl}, \Cref{thm:bound-subfield}, and \Cref{thm:containment} give certain constraints on the successive minima of orders in number fields. We now show that in the limit, these are \emph{all} the constraints. 

\begin{definition}
\label{def:po}
To an order $\cO$ in a degree $n$ number field, associate the point 
\[
	p_{\cO} \coloneqq (\log_{ \lvert \Disc(\cO) \rvert }\lambda_{1}(\cO),\dots,\log_{ \lvert \Disc(\cO) \rvert }\lambda_{n-1}(\cO)) \in \RR^{n-1}.
\]
\end{definition}

\begin{definition}
\label{def:spectrum}
Given a set $\Sigma$ of orders in degree $n$ number fields, let $\Spe(\Sigma)$ denote the set of limit points of the multiset $\{p_{\cO}\}_{\cO \in \Sigma}$.
\end{definition}

\noindent Observe that
\begin{equation*}
\label{eqn:basic-constraint}
\Spe(\Sigma) \subseteq \{\mathbf{x} \in \RR^{n-1} : \sum_{i=1}^{n-1}x_i = 1/2 \text{ and } 0 \leq x_1 \leq \dots \leq x_{n-1}\}.
\end{equation*}
This assertion follows from Minkowski's second theorem, which implies that $\prod_{i = 1}^{n-1} \lambda_i(\cO) \asymp_n \lvert \Disc(\cO) \rvert^{1/2}$, and the fact that $1 \leq \lambda_1 \leq \dots \leq \lambda_{n-1}$.

\begin{definition}
\label{def:sigma-g}
For a permutation group $G \subseteq S_n$, let $\Sigma(G)$ denote the set of (isomorphism classes of) orders in degree $n$ number fields with Galois group $G$. Let $\Sigma_n$ denote the set of (isomorphism classes of) orders in degree $n$ number fields.     
\end{definition}

We would like to compute $\Spe(\Sigma(G))$ and $\Spe(\Sigma_n)$. Our previous theorems (\Cref{thm:bl}, \Cref{thm:bound-subfield}, and \Cref{thm:containment}) imply that $\Spe(\Sigma(G))$ and $\Spe(\Sigma_n)$ are contained in certain linear half-spaces. For example, letting $x_1,\dots,x_{n-1}$ be the coordinates of $\RR^{n-1}$, \Cref{thm:bl} implies that $\Spe(\Sigma(S_n))$ is contained in the linear half-space $x_{i+j} \leq x_i + x_j$ for all $1 \leq i,j < i+j < n$. Our next theorem shows that $\Spe(\Sigma(S_n))$ is (essentially) \emph{equal} to the intersection of these linear half-spaces. 

\begin{theorem}
\label{thm:sn-spectrum}
$\Spe(\Sigma(S_n))$ consists of the points $(x_1,\dots,x_{n-1}) \in \RR^{n-1}$ such that:
\begin{enumerate}
    \item $\sum_{i = 1}^{n-1}x_i = 1/2$;
    \item $0 \leq x_1 \leq x_2 \leq \dots \leq x_{n-1}$;
    \item and $x_{i+j} \leq x_i + x_j$ for all $1 \leq i,j < i+j < n$.
\end{enumerate}
\end{theorem}

In general, we prove that (\Cref{thm:mult-spectrum}) $\Spe(\Sigma_n)$ is a finite union of polytopes. (In this paper, a polytope is the intersection of finitely many linear half-spaces). Lenstra conjectured (\Cref{conj:lenstra}) an explicit description of this finite union of polytopes; in \Cref{thm:lenstra-spectrum}, we'll show that when $n$ is a prime power, a product of $2$ primes, or $12$, Lenstra's conjecture is true. For all other $n$, Lenstra's conjecture is false. To state Lenstra's conjecture, we first introduce some notation.

\begin{definition}
\label{def:lenstra-polytope}
The \emph{Lenstra polytope} $\Len_{\mfT}$ of a tower type $\mfT$ is the set of $\mathbf{x} = (x_1,\dots,x_{n-1}) \in \RR^{n-1}$ satisfying the following conditions:
\begin{enumerate}
\item $\sum_{i = 1}^{n-1}x_i = 1/2$;
\item $0 \leq x_1 \leq x_2 \leq \dots \leq x_{n-1}$;
\item and $x_{i+j} \leq x_i + x_j$ for $i+j$ not overflowing modulo $\mfT$.
\end{enumerate}
\end{definition}

\begin{conjecture}[Lenstra]
\label{conj:lenstra}
\[
	 \Spe(\Sigma_n) = \bigcup_{\mfT} \Len_{\mfT}.
\]
\end{conjecture}

\begin{theorem}
\label{thm:lenstra-spectrum}
If $n$ is a prime power, a product of two primes, or $12$, then
\[
    \Spe(\Sigma_n) = \bigcup_{\mfT} \Len_{\mfT}
\]
If $n$ is not a prime power, a product of two primes, or $12$, then $\Spe(\Sigma_n)$ strictly contains $\cup_{\mfT} \Len_{\mfT}$.
\end{theorem}

Note that $\Spe(\Sigma_n)$ is not always convex! For example, when $n = 6$, the region $\Spe(\Sigma_n)$ is a union of two polytopes. 

\begin{question}
If $n$ is not a prime power, a product of two primes, or $12$, then what is $\Spe(\Sigma_n)$?
\end{question}

We now state a general theorem which shows that $\Spe(\Sigma(G))$ is a finite union of polytopes, beginning with some notation. Let $K/L$ be a degree $n$ field extension.

\begin{definition}
\label{def:flag}
A \emph{flag of $K/L$} is a set $\cF = \{F_0,\dots,F_{n-1}\}$ of $L$-vector spaces such that $L = F_0 \subset F_1 \subset \dots \subset F_{n-1} = K$ and $\dim_{L} F_i = i + 1$ for all $i \in [n]$.  
\end{definition}

\begin{definition}
\label{def:flag-type}
A flag type is a function $T \colon [n]\times[n] \rightarrow [n]$ such that:
\begin{enumerate}
    \item $T(i,j) = T(j,i)$ for all $i,j \in [n]$;
    \item $T(0,i) = i$ for all $i \in [n]$;
    \item and $T(i-1,j) \leq T(i,j)$ for all $j \in [n]$ and all $1 \leq i < n$.
\end{enumerate}
\end{definition}

\begin{definition}
\label{def:tf}
To a flag $\cF$, associate the flag type $T_{\cF}$ given by the formula: 
\begin{align*}
  T_{\cF} \colon [n]\times[n] & \longrightarrow [n] \\
  (i\;\;,\;\;j ) \; & \longmapsto \min\{k \in [n] : F_iF_j \subseteq F_k\}.
\end{align*}
\end{definition}

\begin{definition}
\label{def:ptf}
Given a flag type $T \colon [n]\times[n] \rightarrow [n]$, the polytope $P_{T}$ is the set of $\mathbf{x} = (x_1,\dots,x_{n-1}) \in \RR^{n-1}$ satisfying the following conditions:
\begin{enumerate}
\item $\sum_{i=1}^{n-1} x_i = 1/2$;
\item $0 \leq x_1 \leq \dots \leq x_{n-1}$;
\item and $x_{T(i,j)} \leq x_i + x_j$ for $1 \leq i,j < n$.
\end{enumerate}
\end{definition}

\begin{theorem}
\label{thm:mult-spectrum}
We have
\[
	\Spe(\Sigma(G)) = \bigcup_{\cF}P_{T_{\cF}}
\]
where $\cF$ ranges across all flags of degree $n$ number fields with Galois group $G$.
\end{theorem}

The proofs of \Cref{thm:sn-spectrum} and \Cref{thm:lenstra-spectrum} involve computing $\bigcup_{\cF}P_{T_{\cF}}$ and then applying \Cref{thm:mult-spectrum}.

\subsection{Bounds on scrollar invariants of curves}
We now switch focus and discuss scrollar invariants of curves. Let $k$ be a field and let $C$ be a smooth projective geometrically irreducible curve over $k$ equipped with a finite morphism $\pi \colon C \rightarrow \PP^1$ of degree $n$. Let $\cL$ be a line bundle on $C$. 

\begin{definition}
Let $e_0(\cL) \leq e_1(\cL) \leq \dots \leq e_{n-1}(\cL)$ be the unique integers such that 
\[
	\pi_*\cL \simeq \cO_{\PP^1}(-e_0(\cL)) \oplus \cO_{\PP^1}(-e_1(\cL)) \oplus \dots \oplus \cO_{\PP^1}(-e_{n-1}(\cL)).
\]
We say $e_i$ is the $i$th scrollar invariant of $\cL$ with respect to $\pi$. 
\end{definition}

\begin{theorem}
\label{thm:bl-fn-field}
If $\pi$ doesn't factor through any nontrivial proper subcovers, then for any three line bundles $\cL_1,\cL_2,\cL_3$ with $\cL_1 \otimes \cL_2 \simeq \cL_3$, we have
\[
    e_{i + j}(\cL_3) \leq e_i(\cL_1) + e_i(\cL_2)
\]
for any integers $0 \leq i,j \leq i + j < n$.
\end{theorem}

\begin{theorem}
\label{thm:bound-subfield-fn-field}
Choose integers $0 \leq i,j \leq i + j < n$. Suppose that for every integer $m$ such that $\pi$ factors through a degree $m$ subcover, we have $(i\%m) + (j\%m) = (i + j)\%m$. Then for any three line bundles $\cL_1,\cL_2,\cL_3 \in \Pic(C)$ with $\cL_1 \otimes \cL_2 \simeq \cL_3$, we have
\[
   e_{i + j}(\cL_3) \leq e_i(\cL_1)+e_j(\cL_2).
\]
\end{theorem}

\begin{theorem}
\label{thm:containment-fn-field}
Suppose $n$ is a prime power, a product of $2$ primes, or equal to $12$.  Then there exists a tower type $\mfT$, dependent only on $\pi$, such that for all $0 \leq i,j \leq i + j < n$, if $i + j$ does not overflow modulo $\mfT$, then
\[
    e_{i + j}(\cO_C) \leq e_i(\cO_C) + e_j(\cO_C).
\]
\end{theorem}

\subsection{Previous work}
In the case of successive minima, our results are inspired by and generalize work of Chiche-lapierre, who computed the successive minima spectrum, in different language, for $n = 3,4$ \cite{chiche}; work of Bhargava, Shankar, Taniguchi, Thorne, Tsimerman, and Zhao, and independently Pikert and Rosen, who proved that $\lambda_{n-1} \ll \lambda_i\lambda_j$ for all $i + j = n-1$ \cite{TwoTorsion, pikert}; and unpublished work of Bhargava and Lenstra, who proved that $\lambda_{i+j} \ll \lambda_i\lambda_j$ for all $1 \leq i \leq j \leq i+j < n$ for orders in primitive number fields. In the case of scrollar invariants, our work generalizes classical bounds on the Maroni invariant of trigonal covers \cite{maroni}; results of Ohbuchi bounding the sum of scrollar invariants \cite{obuchi}; and results of Deopurkar and Patel bounding the smallest scrollar invariant \cite{anands}. Combined with recent work of Castryk, Vermeulen, and Zhao \cite{wouter}, our work also provides new constraints on the syzygy bundles of curves.

Related questions have been also addressed by Terr \cite{Terr97}, who proved the equidistribution of shapes of cubic fields; by Bhargava and H \cite{BhaHa16}, who proved the equidistribution of shapes of $S_n$-fields for $n = 4,5$; and by Holmes \cite{holmes}, who proved the equidistribution of shapes in pure prime degree number fields. Our approach differs from that of Terr, Bhargava, H, and Holmes in the following meaningful sense; for $n = 3,4,5$, when ordered by absolute discriminant, the theorems of Bhargava and H imply that $100\%$ of orders in $S_n$-fields lie ``near'' the point $\frac{1}{2(n-1)}(1,\dots,1)$. Thus, equidistribution theorems only ``see'' that one point, but give very refined information at that point. Conversely, our work is focused on classifying the full spectrum of successive minima that may occur, even if much of the spectrum occurs with density $0$.

\subsection{Outline}
In \Cref{sec:constraints}, we introduce a theorem from additive combinatorics and use it to prove bounds on successive minima. Along the way, we provide a proof of \Cref{thm:bound-subfield} and \Cref{prop:mult-bound}. In \Cref{sec:joint-constraints}, we build upon the aforementioned theorem from additive combinatorics to prove joint constraints on successive minima. In particular we prove \Cref{thm:containment}. In \Cref{sec:construction}, we give a construction of orders with almost prescribed successive minima, showing that the constraints arising in \Cref{sec:joint-constraints} are ``all'' the constraints. This construction, along with the work in \Cref{sec:joint-constraints}, gives a proof of \Cref{thm:sn-spectrum} and \Cref{thm:mult-spectrum}. Next, in \Cref{sec:lenstra-spectrum-good}, we explicitly compute the successive minima spectrum when $n$ is a prime power, a product of $2$ primes, or $12$. To continue, in \Cref{sec:lenstra-spectrum-bad}, we explicitly show the successive minima spectrum is larger than conjectured in \Cref{conj:lenstra} when $n$ is not a prime power, a product of $2$ primes, or $12$. Combined with the previous section, this gives a proof of \Cref{thm:lenstra-spectrum}. Finally, in \Cref{sec:fn-field}, we prove bounds on scrollar invariants of curves using the tools built in \Cref{sec:constraints} and \Cref{sec:joint-constraints}. Namely, we prove \Cref{thm:bound-subfield-fn-field} and \Cref{thm:containment-fn-field}.

\subsection{Acknowledgments}
I am extremely grateful to Hendrik Lenstra for the many invaluable ideas, conversations, and corrections throughout the course of this project. I also thank Manjul Bhargava for suggesting the questions that led to this paper and for providing invaluable advice and encouragement throughout the course of this research. Thank you as well to Jacob Tsimerman, Akshay Venkatesh, and Arul Shankar for feedback and illuminating conversations. The author was supported by the NSF Graduate Research Fellowship. 

\section{Constraints on successive minima}
\label{sec:constraints}

The goal of this section is to provide a proof of \Cref{thm:bound-subfield} and \Cref{prop:mult-bound}. Along the way, we'll introduce one of the main technical inputs in this article (\Cref{lem:overflow}). We begin with two elementary lemmas on successive minima. 

\begin{lemma}
\label{lem:lambda-zero}
If $\cO$ is an order in a number field, then $\lambda_0(\cO) = 1$.
\end{lemma}
\begin{proof}
Suppose $\cO$ is an order in a number field $K$ of degree $n$. Note that $\lvert 1 \rvert = 1$ so $\lambda_0 \leq 1$. Let $\sigma_1,\dots,\sigma_{n}$ be the nonzero homomorphisms of $\cO$ into the complex numbers. Then for any nonzero $v \in \cO$,
\begin{align*}
\lvert v \rvert^2 &= \frac{1}{n}\bigg(\sum_{i = 1}^{n} \lvert \sigma_i(v)\rvert^2\bigg) \\
&\geq \sqrt[n]{\prod_{i = 1}^{n} \lvert \sigma_i(v)\rvert^2} && \text{ by the AM-GM inequality} \\
&= \sqrt[n]{\prod_{i = 1}^{n}  \sigma_i(v)^2} \\
&= \lvert N_{K/\QQ}(v)\rvert^{2/n}  \\
& \geq 1.
\end{align*}
Thus, $\lambda_0 \geq 1$.
\end{proof}

\begin{lemma}
\label{lem:sqrt-n}
Suppose we have $u,v \in K$ for some number field $K$ of degree $n$. Then $\lv uv \rv \leq \sqrt{n}\lv u \rv \lv v \rv$.
\end{lemma}
\begin{proof}
We have
\begin{align*}
\lv uv \rv^2 &= \frac{1}{n}\sum_{i = 1}^n\lv \sigma_i(uv)\rv^2 \\
&= \frac{1}{n}\sum_{i = 1}^n\lv \sigma_i(u)\rv^2\lv \sigma_i(v)\rv^2 \\
&\leq \frac{1}{n}\Big (\sum_{i = 1}^n\lv \sigma_i(u)\rv^2 \Big)\Big(\sum_{i = 1}^n\lv \sigma_i(v)\rv^2\Big) \\
&= n\Big (\frac{1}{n}\sum_{i = 1}^n\lv \sigma_i(u)\rv^2 \Big)\Big(\frac{1}{n}\sum_{i = 1}^n\lv \sigma_i(v)\rv^2\Big) \\
&= n\lv u\rv^2\lv v\rv^2.
\end{align*}
\end{proof}

\begin{proof}[Proof of \Cref{prop:mult-bound}]
Let $i,j,k,$ be as in the statement of \Cref{prop:mult-bound}. Set
\[
    S \coloneqq \{u_{i'}v_{j'} : i' \leq i, j' \leq j\}.
\]
By assumption, the vectors in $S$ span a vector space of dimension $k+1$ and are contained in the fractional ideal $\mfa_3$. Because we have exhibited at least $k+1$ linearly independent elements of $\mfa_3$, we have:
\begin{align*}
\lambda_k(\mfa_3) &\leq \max\{\lv v \rv : v \in S\} && \\
&\leq \max\{\sqrt{n}\lv v_{i'} \rv \lv u_{j'} \rv : i' \leq i, j' \leq j\} && \text{ by \Cref{lem:sqrt-n}} \\
&= \max\{\sqrt{n}\lambda_{i'}(\mfa_1)\lambda_{j'}(\mfa_2) \rv : i' \leq i, j' \leq j\} && \\
&= \sqrt{n}\lambda_{i}(\mfa_1)\lambda_{j}(\mfa_2). && \\
\end{align*}
\end{proof}

We will need the following theorem. Given a field extension $K/L$ and two $L$-vector spaces $I,J \subseteq K$, let $IJ$ denote the $L$-vector space $\{vu : v \in I, \; u \in J\}$, where multiplication is multiplication in the field $K$. Let $\Stab(IJ) \coloneqq \{v \in K : vIJ = IJ\}$, where the action of $v$ on $IJ$ is multiplication in $K$.

\begin{theorem}[Bachoc, Serra, Z\'emor \cite{zemorthm}, Theorem~3]
\label{thm:zemor}
Let $K/L$ be a field extension and let $I \subseteq K$ be a finite-dimensional $L$-vector space. There exists a subfield $F_I \subseteq K$ with $F_I \neq L$ such that for each finite-dimensional $L$-vector space $J \subseteq K$, precisely one of the following happens:
\begin{enumerate}
\item $\dim_L IJ \geq \dim_L I + \dim_L J  - 1$;
\item or $\dim_L IJ < \dim_L I + \dim_L J  - 1$ and $F_I IJ = IJ$.
\end{enumerate}
\end{theorem}

\noindent We will use the following corollary of \Cref{thm:zemor}. 

\begin{corollary}
\label{lem:overflow}
Let $K/L$ be a field extension of degree $n$. Choose integers $0 \leq i,j,i+j < n$. Let $I,J$ be dimension $i+1$ (resp. $j+1$) $L$-vector spaces in $K$ and suppose $\dim_{L}IJ \leq i+j$. Set $F \coloneqq \Stab(IJ)$ and $m \coloneqq [F : L]$ and write $i$ and $j$ in mixed radix notation with respect to $(m, n/m)$ as
\begin{align*}
i &= i_1 + i_2 m \\
j &= j_1 + j_2 m.
\end{align*}

Then $m > 1$, $i_1 + j_1 \geq m$, $\dim_{L} FI = (i_2 + 1)m$, $\dim_{L} FJ = (j_2 + 1)m$, and $\dim_{L} IJ = (i_2 + j_2 + 1)m$.
\end{corollary}
\begin{proof}
By assumption, $\dim_{L}IJ \leq i+j = \dim_L I + \dim_L - 1$. So in the notation of \Cref{thm:zemor}, we have $F_I IJ = IJ$, so $F_I \subseteq F$. Therefore $F$ is nontrivial, so $m > 1$. 

\noindent \textbf{Case 1: $F = K$.}
If $F = K$, then because $FIJ = IJ$, we have $IJ = K$, so $\dim_L IJ = n$. By assumption $\dim_{L}IJ \leq i + j < n$, which is a contradiction.

\noindent \textbf{Case 2: $F \neq K$.} 
First, we will need the following claim.

\noindent \textbf{Claim: $\dim_F IJ \geq  \dim_F FI + \dim_F FJ - 1$.}
Assume for the sake of contradiction that the claim is false. Then \Cref{thm:zemor}, applied to the extension $K/F$, implies that there exists a field $M$ strictly containing $F$ such that $M \subseteq \Stab(FIFJ)$. Because $F = \Stab(IJ)$, we have $FIFJ = IJ$; hence $M \subseteq \Stab(IJ) = F$, which is a contradiction. 

Proceeding with the proof, we have:
\begin{align}
i_1 + i_2m + j_1 + j_2m + 1 &= i + j + 1 && \\
&> \dim_L IJ && \\
&= (\dim_L F)(\dim_F IJ) && \text{because $IJ$ is an $F$-vector space} \\
&= m(\dim_F IJ) &&  \\
& \geq m(\dim_F FI) + m(\dim_F FJ) - m && \text{by the claim}\\
& \geq m\bigg \lceil \frac{i+1}{m}\bigg \rceil + m\bigg \lceil \frac{j+1}{m}\bigg \rceil - m && \\
&= m(i_2 + 1) + m(j_2 + 1) - m && \\
&= i_2m + j_2m + m. &&
\end{align}
The inequality $i_1 + i_2m + j_1 + j_2m + 1 > i_2m + j_2m + m$ implies that $i_1 + j_1 \geq m$. 

By definition, $i_1,j_1 < m$, so $i_1 + j_1 < 2m$. Therefore, after rounding $i_1 + i_2m + j_1 + j_2m + 1$ up to the nearest $m$th multiple, we get $m(i_2 + j_2 + 1)$, but this is precisely line $(8)$ of the inequality above. Hence, 
\[
    \dim_L IJ = i_2m + j_2m + m,
\]
and the inequalities from line $(3)$ to line $(8)$ are all equalities.

In particular, the inequality on line $(6)$ of the calculation above must be an equality, so:
\[
	\dim_F FI = \bigg \lceil \frac{i+1}{m}\bigg \rceil
\]
\[
	\dim_F FJ = \bigg \lceil \frac{j+1}{m}\bigg \rceil.
\]
\end{proof}

\begin{corollary}
\label{cor:zemor-coro}
With the notation of \Cref{lem:overflow}, we have
\[
    \dim_{L}IJ \geq \dim_L I + \dim_L J - \dim_L(\Stab(IJ)).
\]
If $\dim_{L}IJ < \dim_L I + \dim_L J  - 1$, then $i + j$ overflows modulo $\dim_L(\Stab(IJ))$.
\end{corollary}
\begin{proof}
If
\[
    \dim_{L}IJ \geq \dim_L I + \dim_L J - 1,
\]
then the assertion is trivially true, as $\dim_L(\Stab(IJ)) \geq 1$. If 
\[
     \dim_{L}IJ < \dim_L I + \dim_L J - \dim_L(\Stab(IJ))
\]
then in the notation of \Cref{lem:overflow}
\begin{align*}
\dim_L IJ &= i_2m + j_2m + m && \text{by \Cref{lem:overflow}}\\
&\geq i_2m + j_2m + m - (2m - i_1 - j_1 - 2) && \text{because $i_1,j_1 < m$} \\
&= (i+1) + (j+1) - m && \\
&= \dim_L I + \dim_L J - m.
\end{align*}
\end{proof}

\begin{definition}
For a fractional ideal $\mfa$, we say $\{v_0,\dots,v_{n-1}\} \subseteq \mfa$ is a set of successive minima representatives for $\mfa$ if the $v_i$ are linearly independent and $\lvert v_i \rvert = \lambda_i(\mfa)$ for all $i \in [n]$.
\end{definition}

\begin{proof}[Proof of \Cref{thm:bound-subfield}]
Let $v_0,\dots,v_{n-1}$ (resp. $u_0,\dots,u_{n-1}$) be successive minima representatives for $\mfa_1$ (resp. $\mfa_2$). Set $I \coloneqq \QQ\langle v_0,\dots,v_i\rangle$ and $J \coloneqq \QQ\langle u_0,\dots,u_j\rangle$. If $\dim_{\QQ}IJ \geq i + j + 1$, then \Cref{prop:mult-bound} implies that
\[
    \lambda_{i+j}(\mfa_3) \leq \lambda_i(\mfa_1)\lambda_j(\mfa_2),
\]
which is the desired conclusion. 

Now assume for the sake of contradiction that $\dim_{\QQ}IJ \leq i + j$ and set $m = \dim_{\QQ}\Stab(IJ)$.  The conclusion of \Cref{lem:overflow} states that $i_1 + j_1 \geq m$. Therefore,   
\[
    (i\%m) + (j\%m) \neq (i + j)\%m.
\]
However, this contradicts the assumptions of \Cref{thm:bound-subfield} because $\Stab(IJ)$ is a field.
\end{proof}

\section{Joint constraints on successive minima}
\label{sec:joint-constraints}

As we've shown in \Cref{sec:constraints}, multiplication induces constraints on the successive minima of fractional ideals in number fields. It is natural to ask: how do these constraints interact with each other? In this section we address this question by providing a proof of \Cref{thm:containment}.

The key observation on joint constraints on successive minima is the following. Let $\cO$ be an order in a degree $n$ number field. It is known (see, e.g., \cite{siegel}, Lecture 10, \S 6) that there exists a \emph{Minkowski reduced basis} $\{v_0=1, v_1, v_2,\dots,v_{n-1}\}$ for $\cO$ such that
\begin{equation}
\label{eqn:red}
    \lambda_i(\cO) \asymp_n \lvert v_i \rvert
\end{equation}
and for every $v = \sum_{i = 0}^{n-1}c_i v_i \in \cO$, we have
\begin{equation}
\label{eqn:size}
    \lvert v\rvert \asymp_n \sum_{i = 0}^{n-1}\lvert c_i \rvert \lambda_i(\cO).
\end{equation}

Let $\cF = \{F_i\}_{i \in [n]}$ be the corresponding flag; that is, let $F_i \coloneqq \QQ \langle 1=v_0,v_1,\dots,v_i\rangle$. Let $T_{\cF}$ be the flag type (see \Cref{def:tf}) corresponding to $\cF$.

\begin{proposition}
\label{prop:flag-bound}
For every $0 \leq i,j < n$, we have
\[
    \lambda_{T_{\cF}(i,j)}(\cO) \ll_n \lambda_{i}(\cO)\lambda_j(\cO).
\]
\end{proposition}
\begin{proof}
Let $k = T_{\cF}(i,j)$. By definition, $k$ is the smallest integer such that $F_i F_j \subseteq F_k$. The vector space $F_iF_j$ is spanned by the set
\[
    S \coloneqq \{v_{i'}v_{j'} : i' \leq i, j' \leq j\},
\]
so there exists some $i'\leq i$ and $j' \leq j$ such that the basis expansion
\[
    v_{i'}v_{j'} = \sum_{i = 1}^{n-1}c_i v_i
\]
has $c_k \neq 0$. Therefore we have: 
\begin{align*}
\lambda_k(\cO) &\ll_n \lvert v_{i'}v_{j'} \rvert && \text{because $c_k \neq 0$ and \Cref{eqn:size}} \\
&\ll_n \lvert v_{i'} \rvert \lv v_{j'} \rv && \text{by \Cref{lem:sqrt-n}} \\
&\asymp_n \lambda_{i'}(\cO)\lambda_{j'}(\cO) && \text{by \Cref{eqn:red}} \\
& \leq \lambda_{i}(\cO)\lambda_{j}(\cO).
\end{align*}
\end{proof}

So, to understand joint constraints on successive minima, it is necessary to understand the combinatorics of the flag types $T_{\cF}$. Towards this goal, our main technical result is \Cref{thm:flag-types}, which we prove in \Cref{subsec:flag-types}. To state this theorem, we first introduce some notation. 

\begin{definition}
\label{def:T-mft}
Fix a tower type $\mfT = (n_1,\dots,n_t)$. For any integers $0 \leq i , j < n$, write
\[
	i = i_1 + i_2 n_1 + i_3 (n_1n_2) + \dots + i_t(n_1\dots n_{t-1})
\]
\[
	j = j_1 + j_2 n_1 + i_3 (n_1n_2) + \dots + j_t(n_1\dots n_{t-1})
\]
in mixed radix notation with respect to $\mfT$. For $1 \leq \ell \leq t$, set $k_{\ell} \coloneqq \min(n_{\ell}-1, i_{\ell} + j_{\ell})$. Define the tower type $T_{\mfT}$ by
\[
    T_{\mfT}(i,j) \coloneqq k_1 + k_2 n_1 + k_3 (n_1n_2) + \dots + k_t(n_1\dots n_{t-1}).
\]
\end{definition}

It is easy to see that $T_{\mfT}$ is a flag type: it trivially satisfies properties $(1)$ and $(2)$ of \Cref{def:flag-type}, and an easy calculation shows that $T_{\mfT}$ satisfies property $(3)$ as well. Next, we endow the set of flag types with a poset structure. 

\begin{definition}
For any two flag types $T$ and $T'$, say $T \leq T'$ if $T(i,j) \leq T'(i,j)$ for all $i,j \in [n]$.
\end{definition}

\begin{theorem}
\label{thm:flag-types}
If $n$ is a prime power, a product of two primes, or $12$, then for every flag $\cF$ of a degree $n$ field extension, there exists a tower type $\mfT$ such that $T_{\mfT} \leq T_{\cF}$.
\end{theorem}
The proof of \Cref{thm:flag-types} can be found in \Cref{subsec:flag-types}. Finally, to use \Cref{thm:flag-types}, we need to understand the flag types $T_{\mfT}$.
\begin{lemma}
\label{cor:explicit-flag-type-basic}
For any tower type $\mfT$ and $1 \leq i,j,i+j < n$, the following are equivalent:
\begin{enumerate}
\item $i + j$ does not overflow modulo $\mfT$;
\item and $T_{\mfT}(i,j) = i+j$.
\end{enumerate}
\end{lemma}

 We delay the proof of \Cref{cor:explicit-flag-type-basic} to \Cref{subsec:explicit-flag-type}, where we prove a generalization (\Cref{cor:explicit-flag-type}). Now, we can finally provide a proof of \Cref{thm:containment}, assuming \Cref{thm:flag-types} and \Cref{cor:explicit-flag-type-basic}.

\begin{proof}[Proof of \cref{thm:containment}]
Let $\cF$ be a flag obtained from a Minkowski reduced basis of $\cO$. By \Cref{thm:flag-types}, there exists a tower type $\mfT$ such that $T_{\mfT} \leq T_{\cF}$. Let $0 \leq i,j \leq i+j < n$ be integers such that $i+j$ does not overflow modulo $\mfT$. Then:
\begin{align*}
\lambda_{i+j} &=  \lambda_{T_{\mfT}(i,j)} && \text{ because $i+j = T_{\mfT}(i,j)$ by  \Cref{cor:explicit-flag-type-basic}} \\
&\leq \lambda_{T_{\cF}(i,j)} && \text{because $T_{\mfT} \leq T_{\cF}$, so $T_{\mfT}(i,j) \leq T_{\cF}(i,j)$}\\
&\ll_n \lambda_i \lambda_j && \text{by \Cref{prop:flag-bound}}
\end{align*}
\end{proof}

\subsection{Explicit description of the flag type $T_{\mfT}$}
\label{subsec:explicit-flag-type}

The purpose of this subsection is to explicitly describe the flag types $T_{\mfT}$ by proving \Cref{cor:explicit-flag-type}, beginning with a crucial definition. 

\begin{definition}
\label{def:corner}
Given a flag type $T$, say $(i,j)$ is a \emph{corner} of $T$ if $0 < i,j < n$ and $T(i-1,j) < T(i,j)$ and $T(i,j-1) < T(i,j)$.
\end{definition}

\begin{lemma}
\label{cor:explicit-flag-type}
For any tower type $\mfT$ and $1 \leq i,j < i+j < n$, the following are equivalent:
\begin{enumerate}
\item $i + j$ does not overflow modulo $\mfT$;
\item $(i,j)$ is a corner of $T_{\mfT}$;
\item and $T_{\mfT}(i,j) = i+j$.
\end{enumerate}
\end{lemma}
\begin{proof}
We first show the equivalence of $(1)$ and $(3)$. Letting $k$ be as in the notation of \Cref{def:T-mft}, we can easily see that $i + j$ does not overflow modulo $\mfT$ if and only if $k_{\ell} = i_{\ell} + j_{\ell}$ for all $\ell$. Now, the definition of $T_{\mfT}$ (\Cref{def:T-mft}) shows that this is equivalent to $T_{\mfT}(i,j) = i+j$.

We now show $(3) \implies (2)$. Choose $1 \leq i,j < n$ so that $T_{\mfT}(i,j) = i+j$; equivalently, $k_{\ell} = i_{\ell} + j_{\ell}$ for all $\ell$. We'll show that $T_{\mfT}(i-1,j) < T_{\mfT}(i,j)$. A completely symmetric argument will show that $T_{\mfT}(i,j-1) < T_{\mfT}(i,j)$. 

Suppose $i_1 \neq 0$. Then
\[
	i - 1 = (i_1-1) + i_2 n_1 + i_3 (n_1n_2) + \dots + i_t(n_1\dots n_{t-1}).
\]
in mixed radix notation. Because $k_{\ell} = i_{\ell} + j_{\ell}$ for all $\ell$, we have (in particular) $k_1 = i_1 + j_1$. Clearly $k_1 \leq n_1 - 1$, so $(i_1 - 1) + j_1 \leq n_1 - 1$. Therefore, by the definition of $T_{\mfT}$, we see that $T_{\mfT}(i-1,j) = i - 1 + j < T_{\mfT}(i,j)$. 

Now suppose $i_1 = 0$ and let $\ell$ be the smallest integer such that $i_{\ell} \neq 0$ (such an integer exists because $i \neq 0$). By assumption $\ell \geq 2$. Then
\[
i - 1 = (n_1-1) + \dots + (n_{\ell-1}-1)(n_1\dots n_{\ell-2}) + (i_{\ell}-1)(n_1\dots n_{\ell-1}) + i_{\ell+1}(n_1\dots n_{\ell}) + \dots + i_{t}(n_1\dots n_{t}).
\]
in mixed radix notation. Therefore, we have
\begin{align*}
T_{\mfT}(i-1,j) &= (n_1-1) + \dots + (n_{\ell-1}-1)(n_1\dots n_{\ell-2}) + (k_{\ell}-1)(n_1\dots n_{\ell-1}) + k_{\ell+1}(n_1\dots n_{\ell}) + \dots + k_{t}(n_1\dots n_{t}) \\
&< k_{\ell}(n_1\dots n_{\ell-1}) + k_{\ell+1}(n_1\dots n_{\ell}) + \dots + k_{t}(n_1\dots n_{t}) \\
&\leq k_1 + k_2 n_1 + k_3 (n_1n_2) + \dots + k_t(n_1\dots n_{t-1}) \\
&= T_{\mfT}(i,j)
\end{align*}

We now show $(2) \implies (3)$ by proving the contrapositive. Choose $1 \leq i,j < n$ so that $T_{\mfT}(i,j) \neq i+j$; then there exists some $\ell$ such that $i_{\ell} + j_{\ell} \geq n_{\ell}$, so $k_{\ell} = n_{\ell} - 1$. Without loss of generality suppose $i_{\ell} \neq 0$; set 
\[
    i' \coloneqq i_1 + i_2 n_1 + \dots + (i_{\ell}-1) (n_1\dots n_{\ell-1}) + \dots + i_t(n_1\dots n_{t-1}).
\]
Because $(i_{\ell}-1) + j_{\ell} \geq n_{\ell} - 1 = k_{\ell}$, we have 
\begin{equation}
\label{eqn:equality}
T_{\mfT}(i',j) = k_1 + k_2 n_1 + \dots + k_t(n_1\dots n_{t-1}) = T_{\mfT}(i,j). 
\end{equation}
By definition,
\[
    T_{\mfT}(i',j) \leq T_{\mfT}(i-1,j) \leq T_{\mfT}(i,j).
\]
so the equality \Cref{eqn:equality} implies that $T_{\mfT}(i-1,j) = T_{\mfT}(i,j)$, so $(i,j)$ is not a corner of $T_{\mfT}$.
\end{proof}

\subsection{Explicit description of the flag types $T_{\cF}$}
\label{subsec:flag-types}
The primary goal of this subsection is to prove \Cref{thm:flag-types}, which is a description of the flag types $T_{\cF}$. We'll first need the following lemma, which we use repeatedly throughout the proof of \Cref{thm:flag-types}.

\begin{lemma}
\label{lem:test-flag-type-inequal}
For any two flag types $T$ and $T'$ such that $T \not \geq T'$, there exists a corner $(i,j)$ of $T'$ such that $T(i,j) < T'(i,j)$.
\end{lemma}
\begin{proof}
Because $T \not \geq T'$, then there exists $(i,j)$ such that $T(i,j) < T'(i,j)$. Choose $i' \leq i$, $j' \leq j$ such that $(i',j')$ is a corner of $T'$ and $T'(i',j') = T'(i,j)$. Because $i' \leq i$ and $j' \leq j$, we have $T(i',j') \leq T(i,j)$. Hence $T(i',j') < T'(i',j')$.
\end{proof}

\begin{proof}[Proof of \Cref{thm:flag-types}]
Follows from combining \Cref{prop:2prefined}, \Cref{prop:pq}, \Cref{prop:prime-power}, and \Cref{prop:12}.
\end{proof}

\begin{proposition}
\label{prop:2prefined}
Suppose $n = 2p$ for $p$ an odd prime. For every flag $\cF$ of a degree $n$ field extension $K/L$, we have $T_{\cF} \geq T_{(2,p)}$ or $T_{\cF} \geq T_{(p,2)}$.
\end{proposition}
\begin{proof}
Assume for the sake of contradiction that there exists a flag $\cF$ of $K/L$ such that $T_{\cF} \not \geq T_{(p,2)}$ and $T_{\cF} \not \geq T_{(2,p)}$. \Cref{lem:test-flag-type-inequal} implies that there exist integers $0 < i_2 \leq j_2 < i_2 + j_2 < 2p$ such that $T_{\cF}(i_2,j_2) < T_{(p,2)}$ and $(i_2,j_2)$ is a corner of $T_{(p,2)}$. Because  $(i_2,j_2)$ is a corner of $T_{(p,2)}$, \Cref{cor:explicit-flag-type} implies that $T_{(p,2)}(i_2,j_2) = i_2+j_2$. Therefore, we have
\begin{equation}
\label{eqn:2p-refined-1}
T_{\cF}(i_2,j_2) < i_2 + j_2.
\end{equation}
Recall that by definition, $T_{\cF}(i,j)$ is equal to the smallest value of $k$ such that $F_{i}F_j \subseteq F_k$, where $F_i$ are the vector subspaces comprising the flag $\cF$. So \Cref{eqn:2p-refined-1} implies that $F_{i_2}F_{j_2} \subseteq F_{i_2 + j_2 - 1}$, so 
\[
	\dim_{L}F_{i_2}F_{j_2} \leq i_2 + j_2.
\]
Now, \Cref{cor:zemor-coro} implies that $i_2 + j_2$ overflows modulo $\dim_{L}\Stab(F_{i_2}F_{j_2})$. Now, because  $(i_2,j_2)$ is a corner of $T(p,2)$, \Cref{cor:explicit-flag-type} implies that $i_2 + j_2$ does not overflow modulo $p$. Because $\Stab(F_{i_2}F_{j_2})$ is a field, its degree over $L$ must be $1$, $2$, $p$, or $2p$; because $i_2 + j_2$ overflows modulo the degree, we have that $\dim_{L}\Stab(F_{i_2}F_{j_2}) = 2$. 

Furthermore, because $i_2 + j_2$ must overflow modulo $2$, \Cref{lem:overflow} implies that
\begin{equation}
\label{eqn:2p-refined-1-mod}
\dim_{L}\Stab(F_{i_2}F_{j_2})F_{i_2} = (\dim_{L}F_{i_2} - 1)2 = \dim_{L}F_{i_2}.
\end{equation}
Hence $F_{i_2}$ is a vector space over $\Stab(F_{i_2}F_{j_2})$.

Again \Cref{lem:test-flag-type-inequal} implies that there exist integers $0 < i_p \leq j_p < i_p + j_p < 2p$ such that $T_{\cF}(i_p,j_p) < T(2,p)$ and $(i_p,j_p)$ is a corner of $T(2,p)$. The same reasoning shows that $\Stab(F_{i_p}F_{j_p})$ is a field of degree $p$ over $L$ and $i_p + j_p$ overflows modulo $p$. Because $i_p + j_p$ overflows modulo $p$ and $i_p + j_p < 2p$, we have
\begin{equation}
\label{eqn:2p-refined-2}
\frac{p+1}{2} \leq j_p < p
\end{equation}
so \Cref{lem:overflow} implies that $\dim_{L}\Stab(F_{i_p}F_{j_p})F_{j_p} = p$. Because $1 \in F_{j_p}$, we have
\[
    \Stab(F_{i_p}F_{j_p}) \subseteq \Stab(F_{i_p}F_{j_p})F_{j_p},
\]
and hence $\Stab(F_{i_p}F_{j_p})F_{j_p} = \Stab(F_{i_p}F_{j_p})$. Because $1 \in \Stab(F_{i_p}F_{j_p})$ we have:
\[
F_{j_p} \subseteq \Stab(F_{i_p}F_{j_p})F_{j_p} = \Stab(F_{i_p}F_{j_p}).
\]
Therefore $F_{j_p}$ is contained in the degree $p$ field $\Stab(F_{i_p}F_{j_p})$.

Because $F_{i_2}$ is a vector space over a quadratic field, it cannot be contained in a degree $p$ field. Therefore, $j_p < i_2$ and Therefore, we must have $F_{j_p}  \subset F_{i_2}$. Putting this all together, we obtain:
\begin{align*}
i_2 + 1 &= \dim_L F_{i_2} && \\
&= \dim_L \Stab(F_{i_2}F_{j_2})F_{i_2} && \text{by \Cref{eqn:2p-refined-1-mod}} \\
&\geq \dim_L \Stab(F_{i_2}F_{j_2})F_{j_p} && \text{because  $F_{j_p}  \subseteq F_{i_2}$ } \\
&= 2\dim_{L}F_{j_p} && \text{because $\deg(\Stab(F_{i_2}F_{j_2})) = 2$ and $F_{j_p} \subseteq$ a degree $p$ field}\\
&\geq 2\bigg(\frac{p+1}{2} + 1\bigg)  && \text{by \Cref{eqn:2p-refined-2}}\\
&= p+3. &&
\end{align*}
Therefore, $i_2 \geq p+2$. Now, because $i_2 \leq j_2$, we have that $i_2 + j_2 \geq 2p$, which is a contradiction.
\end{proof}

\begin{proposition}
\label{prop:pq}
Suppose $n = pq$ for $p$ and $q$ distinct odd primes. Then for every flag $\cF$ of a degree $n$ field extension $K/L$, we have $T_{\cF} \geq T_{(p,q)}$ or $T_{\cF} \geq T_{(q,p)}$.
\end{proposition}
\begin{proof}
Assume for the sake of contradiction that there exists a flag $\cF$ of $K/L$ such that $T_{\cF} \not \geq T_{(p,q)}$ and $T_{\cF} \not \geq T_{(q,p)}$. Identically to the proof of \Cref{prop:2prefined}, there exist integers $0 < i_q \leq j_q < i_q + j_q < pq$ such that:
\begin{enumerate}
\item $i_q + j_q$ does not overflow modulo $(p,q)$;
\item $i_q + j_q$ does overflow modulo $(q,p)$;
\item $T_{\cF}(i_q,j_q) < i_q + j_q$ and $\dim_{L}F_{i_q}F_{j_q} \leq i_q + j_q$ and $\dim_{L}\Stab(F_{i_q}F_{j_q}) = q$. 
\end{enumerate} 
Similarly, there exist integers $0 < i_p \leq j_p < i_p + j_p < pq$ such that:
\begin{enumerate}
\item $i_p + j_p$ does not overflow modulo $(q,p)$;
\item $i_p + j_p$ does overflow modulo $(p,q)$;
\item $T_{\cF}(i_p,j_p) < i_p + j_p$ and $\dim_{L}F_{i_p}F_{j_p} \leq i_p + j_p$ and $\dim_{L}\Stab(F_{i_p}F_{j_p}) = p$. 
\end{enumerate} 
Set $K_q = \Stab(F_{i_q}F_{j_q})$ and $K_p = \Stab(F_{i_p}F_{j_p})$. Without loss of generality, suppose $j_q \leq j_p$.

\noindent \textbf{Case 1: $i_q \leq i_p$.} Then:
\begin{align*}
K_q &\subseteq F_{i_q} F_{j_q} && \text{because $K_q = \Stab(F_{i_q} F_{j_q}) \subseteq  F_{i_q} F_{j_q}$}\\
&\subseteq F_{i_p} F_{j_p} && \text{because $i_q \leq i_p$ and $j_q \leq j_p$}
\end{align*}
Now, $F_{i_p}F_{j_p}$ is a $K_p$-vector space. Because $K_q \subseteq F_{i_p} F_{j_p}$, we have $K_pK_q \subseteq F_{i_p} F_{j_p}$, which implies $K = F_{i_p} F_{j_p}$. But then $\dim_{L}\Stab(F_{i_p} F_{j_p}) = pq$, contradiction. 

\noindent \textbf{Case 2: $i_p \leq i_q$.} Write $i_p$ and $j_p$ in mixed radix notation with respect to $(p,q)$ and write $i_q$ and $j_q$ in mixed radix notation with respect to $(q,p)$ as
\begin{align*}
i_p &= i_{1,p}p + i_{2,p} \\
j_p &= j_{1,p}p + j_{2,p} \\
i_q &= i_{1,q}q + i_{2,q} \\
j_q &= j_{1,q}q + j_{2,q}.
\end{align*} 
By \cref{lem:overflow}, we have that

\begin{align}
\label{eqn:dim-comp-pq-jp}
\dim_L K_p F_{j_p} &= (\dim_L K_p)(j_{1,p}+1) \\
\label{eqn:dim-comp-pq-iq}
\dim_L K_q F_{i_q} &= (\dim_L K_q)(i_{1,q}+1) \\
\label{eqn:dim-comp-pq-jq}
\dim_L K_q F_{j_q} &= (\dim_L K_q)(j_{1,q}+1).
\end{align}

We'll need the following lemma.
\begin{lemma}
\label{lem:helper}
Given a field extension $K/L$, an $L$-vector space $V \subseteq K$, and two subfields $M_1,M_2$ with $M_1 \cap M_2 = L$ and $M_1M_2 = K$, we have that: 
\[
	\dim_{L} V \leq \frac{\dim_{L}VM_1}{\dim_{L}M_1}\frac{\dim_{L}VM_2}{\dim_{L}M_2}. 
\]
\end{lemma}
\begin{proof}
Let $\{\alpha_1,\dots,\alpha_r\}$ be an $M_1$-basis for $VM_1$, and let $\{\beta_1,\dots,\beta_s\}$ be an $M_2$-basis for $VM_2$. 
Extend so that $\{\alpha_1,\dots,\alpha_t\}$ is an $M_1$-basis for $K$ and $\{\beta_1,\dots,\beta_u\}$ is an $M_2$-basis for $K$.  We claim that the set $\{\alpha_i\beta_j\}_{1 \leq i \leq t, 1 \leq j \leq u}$ is $L$-linearly independent. Indeed, if
\[
	\sum_{i = 1}^t\sum_{j = 1}^u c_{ij}\alpha_i\beta_j = 0
\]
for some $c_{ij} \in L$, then because the $\alpha_i$ are $M_1$-linearly independent, we must have $\sum_{j = 1}^u c_{ij}\beta_j = 0$ for all $i$; now because the $\beta_j$ are $M_2$-linearly independent, we must have $c_{ij} = 0$ for all $i,j$.

Because $M_1M_2 = K$, the $L$-linear span of the set $\{\alpha_i\beta_j\}_{1 \leq i \leq t, 1 \leq j \leq u}$ is equal to $K$, and thus $\{\alpha_i\beta_j\}_{1 \leq i \leq t, 1 \leq j \leq u}$ is an $L$-basis of $K$. Now given $x \in V$, write 
\[
	x = \sum_{i = 1}^t\sum_{j = 1}^u c_{ij}\alpha_i\beta_j.
\]
Because $\{\alpha_1,\dots,\alpha_r\}$ is an $M_1$-basis for $VM_1$, we have $\sum_{j = 1}^u c_{ij}\beta_j = 0$ for all $i > r$; now because the $\beta_j$ are $M_2$-linearly independent, we must have $c_{ij} = 0$ for all $i>r$. Similarly, because $\{\beta_1,\dots,\beta_s\}$ are an $M_2$-basis for $VM_2$, we have $\sum_{i = 1}^t c_{ij}\alpha_i = 0$ for all $j > s$; now because the $\alpha_i$ are $M_1$-linearly independent, we must have $c_{ij} = 0$ for all $j>s$. Therefore, $V$ is contained in the $L$-linear span of $\{\alpha_i\beta_j\}_{1 \leq i \leq r, 1 \leq j \leq s}$, so
\[
	\dim_{L} V \leq rs = (\dim_{M_1} VM_1) (\dim_{M_2} VM_2) =  \frac{\dim_{L}VM_1}{\dim_{L}M_1}\frac{\dim_{L}VM_2}{\dim_{L}M_2}. 
\]

\end{proof}

We now continue with the proof of \Cref{prop:pq}. We obtain:
\begin{gather}
\label{eqn:pq-iq+1}
\begin{aligned}
i_q + 1 &= \dim_L F_{i_q} && \\
&\leq \frac{\dim_L K_q F_{i_q}}{\dim_L K_q}\frac{\dim_L K_p F_{i_q}}{\dim_L K_p } && \text{by \Cref{lem:helper}} \\
&= (i_{1,q}+1)\frac{\dim_L K_p F_{i_q}}{\dim_L K_p } && \text{by \Cref{eqn:dim-comp-pq-iq}} \\
&\leq (i_{1,q}+1)\frac{\dim_L K_p F_{j_p}}{\dim_L K_p } && \text{because $i_q \leq j_q \leq j_p$}\\
&= (i_{1,q}+1)(j_{1,p}+1) && \text{by \Cref{eqn:dim-comp-pq-jp}}
\end{aligned}
\end{gather}
Similarly, we get
\begin{gather}
\label{eqn:pq-jq+1}
\begin{aligned}
j_q + 1 &= \dim_L F_{j_q} && \\
&\leq \dim_L F_{j_q} \\
&\leq \frac{\dim_L K_q F_{j_q}}{\dim_L K_q }\frac{\dim_L K_p F_{j_q}}{\dim_L K_p } && \text{by \Cref{lem:helper}} \\
&= (j_{1,q}+1)\frac{\dim_L K_p F_{j_q}}{\dim_L K_p } &&  \text{by \Cref{eqn:dim-comp-pq-jq}} \\
&\leq (j_{1,q}+1)\frac{\dim_L K_p F_{j_p}}{\dim_L K_p } && \text{because $j_q \leq j_p$} \\
&= (j_{1,q}+1)(j_{1,p}+1) &&  \text{by \Cref{eqn:dim-comp-pq-jp}}\\
\end{aligned}
\end{gather}
Combining \Cref{eqn:pq-iq+1} and \Cref{eqn:pq-jq+1}, we see that:
\begin{gather}
\label{eqn:pq-ineq}
\begin{aligned}
i_q + j_q &< (i_{1,q}+1)(j_{1,p}+1) + (j_{1,q}+1)(j_{1,p}+1) \\
		  &= (i_{1,q}+j_{1,q}+2)(j_{1,p}+1) \\
		  &\leq p(j_{1,p}+1) \\
		  &\leq j_p.
\end{aligned}
\end{gather}
We have:
\begin{align*}
K_q &\subseteq F_{i_q}F_{j_q} && \text{because $K_q = \Stab(F_{i_q}F_{j_q})$ and $1 \in F_{i_q}F_{j_q}$} \\
&\subseteq F_{i_q + j_q - 1} && \text{because $T_{\cF}(i_q,j_q) < i_q + j_q$} \\
&\subseteq F_{j_p} && \text{by \Cref{eqn:pq-ineq}} \\
&\subseteq F_{i_p}F_{j_p} && \text{because $1 \in F_{i_p}$}
\end{align*}
Because $F_{i_p}F_{j_p}$ is a $K_p$-vector space, we have
\[
	K = K_qK_p \subseteq F_{i_p}F_{j_p},
\]
which is a contradiction.
\end{proof}

\begin{proposition}
\label{prop:prime-power}
Suppose $n = p^k$ for $p$ a prime and $k \geq 1$. Then for every flag $\cF$ of a degree $n$ field extension $K/L$, we have $T_{\cF} \geq T_{(p,\dots,p)}$.
\end{proposition}
\begin{proof}
Assume for the sake of contradiction that there exists a flag $\cF$ of a degree $n$ field $K$ such that $T_{\cF} \not \geq T_{(p,\dots,p)}$. By \Cref{lem:test-flag-type-inequal}, as in the proof of \Cref{prop:2prefined}, there exists integers $0 < i \leq j < n$ such that $T_{\cF}(i,j) < T_{(p,\dots,p)}(i,j)$ and $(i,j)$ is a corner of $T_{(p,\dots,p)}$. Because $(i,j)$ is a corner of $T(p,\dots,p)$, \Cref{cor:explicit-flag-type} implies that $i+j$ does not overflow modulo $(p,\dots,p)$ and $T_{(p,\dots,p)}(i,j) = i+j$, and so therefore the addition $i + j$ does not overflow modulo $p^{\ell}$ for any positive integer $\ell$. 

Now, $T(i,j) < i+j$, so $\cF_i\cF_j \subseteq \cF_{i + j - 1}$, implying that 
\[
    \dim_{L}\cF_i \cF_j \leq i + j.
\]
Now by \cref{lem:overflow}, the addition $i + j$ must overflow over some positive integer $m$ such that $m \mid p^k$, which is a contradiction.
\end{proof}

\begin{lemma}
\label{lem:i-leq-m}
Suppose $\cF$ is a flag of a degree $n$ field extension $K/L$ and let $0 < i, j < i + j < n$ be integers such that $\cF_i \cF_j \subseteq \cF_{i+j-1}$. Let $m = \deg_{L}\Stab(F_i F_j)$. If $i < m$, then $F_i \subseteq \Stab(F_iF_j)$.
\end{lemma}
\begin{proof}
Because $\cF_i \cF_j \subseteq \cF_{i+j-1}$, we have
\[
    \dim_{L}\cF_i \cF_j \leq i+j.
\]
Applying \Cref{lem:overflow}, we have that $\dim_{L}\Stab(F_iF_j)F_{i} = m$; because $\Stab(F_iF_j) \subseteq F_{i}$ and $\dim_{L}\Stab(F_iF_j) = m$, we have that $\Stab(F_iF_j)F_{i} = \Stab(F_iF_j)$. Therefore, $F_{i} \subseteq \Stab(F_iF_j)$. 
\end{proof}

\begin{proposition}
\label{prop:12}
Suppose $n = 12$. Then for every flag $\cF$ of a degree $n$ field extension $K/L$, we have $T_{\cF} \geq T_{(3,2,2)}$ or $T_{\cF} \geq T_{(2,3,2)}$ or $T_{\cF} \geq T_{(2,2,3)}$.
\end{proposition}
\begin{proof}
Assume for the sake of contradiction that there exists a flag $\cF$ of a degree $12$ field extension $K/L$ such that $T_{\cF} \not \geq T_{(3,2,2)}$, $T_{\cF} \not \geq T_{(2,3,2)}$, and $T_{\cF} \not \geq T_{(2,2,3)}$. As in the proof of \cref{prop:2prefined}, there exist positive integers $i_1$, $i_2$, $i_3$, $j_1$, $j_2$, $j_3$ such that
\begin{align*}
0 < i_1 \leq j_1 &< i_1 + j_1 < 12 \\
0 < i_2 \leq j_2 &< i_2 + j_2 < 12 \\
0 < i_3 \leq j_3 &< i_3 + j_3 < 12,
\end{align*}
and $i_1 + j_1$ (resp. $i_2 + j_2$, $i_3 + j_3$) does not overflow modulo $(3,2,2)$ (resp. $(2,3,2)$, $(2,2,3)$), and
\begin{align*}
\label{cond}
T_{\cF}(i_1 + j_1) &< i_1 + j_1 \\
T_{\cF}(i_2 + j_2) &< i_2 + j_2 \\
T_{\cF}(i_3 + j_3) &< i_3 + j_3.
\end{align*}

Set:
\[
	K_1 \coloneqq \Stab(F_{i_1}F_{j_1})
\]
\[
	K_2 \coloneqq \Stab(F_{i_2}F_{j_2})
\]
\[
	K_3 \coloneqq \Stab(F_{i_3} F_{j_3}).
\]
and let $m_1 = \dim_L K_1$, let $m_2 = \dim_L K_2$, and let $m_3 = \dim_L K_3 $. Because all $K_{\ell}$ are subfields, we have that $m_{\ell} \mid 12$ for all $\ell = 1,2,3$. By \Cref{lem:overflow}, for all $\ell = 1,2,3$, the addition $i_{\ell} + j_{\ell}$ overflows modulo $m_{\ell}$.

We now enumerate the list of possible triples of positive integers $(i_1,j_1,m_1)$ for which $0 < i_1 \leq j_1 < i_1 + j_1 < 12$, the addition $i_1 + j_1$ overflows modulo $m_1$, and $m_1 \mid 12$, and the addition $i_1 + j_1$ does not overflow modulo $(3,2,2)$. The list is:
\[
	\mathcal{L}_1 = \{(1,1,2), (1,3,2), (1,3,4), (1,7,2), (1,7,4), (1,9,2), (2,3,4), (2,6,4), (3,6,4), (3,7,2), (3,7,4)\}.
\]

Similarly, the list of possible triples $(i_2,j_2,m_2)$ for which $0 < i_2 \leq j_2 < i_2 + j_2 < 12$, the addition $i_2 + j_2$ overflows modulo $m_2$, and $m_2 \mid 12$, and the addition $i_2 + j_2$ does not overflow modulo $(2,3,2)$ is:
\[
	\mathcal{L}_2 = \{(1,2,3), (1,8,3), (2,2,3), (2,2,4), (2,3,4), (2,6,4), (2,7,3), (2,7,4), (2,8,3), (3,6,4)\}.
\]

Finally, the list of possible triples $(i_3,j_3,m_3)$ for which $0 < i_3 \leq j_3 < i_3 + j_3 < 12$, the addition $i_3 + j_3$ overflows modulo $m_3$, and $m_3 \mid 12$, and the addition $i_3 + j_3$ does not overflow modulo $(2,2,3)$ is:
\[
	\mathcal{L}_3 = \{(1,2,3), (1,8,3), (2,4,3), (2,4,6), (2,5,3), (2,5,6), (2,8,3), (3,4,6), (4,4,6), (4,5,3), (4,5,6)\}.
\]

We now show that no combination of integers $(i_1,j_1,m_1)$, $(i_2,j_2,m_2)$, and $(i_3,j_3,m_3)$ from the lists above is possible. Let $(i_1,j_1,m_1)$, $(i_2,j_2,m_2)$, and $(i_3,j_3,m_3)$ be any triples from $\mathcal{L}_1,\mathcal{L}_2,\mathcal{L}_3$ respectively. Choose $v_1 \in K$ such that $F_1 = L\langle 1, v_1\rangle$. 

\noindent\textbf{Claim $(A)$: $\deg(v_1) \mid m_1 \mid 4$.}
Because $i < m$ for every $(i,j,m) \in \cL_1$, we have $i_1 < m_1$. \Cref{lem:i-leq-m} implies that $F_{i_1} \subseteq K_1$, so we see that:
\[
    v_1 \in F_1 \subseteq F_{i_1} \subseteq K_1.
\]
Therefore,  the field $L(v_1)$ is contained in a field of degree $m_1$, so $\deg(v_1) \mid m_1$. Now looking at the list $\cL_1$ shows that $m_1 \in \{2,4\}$, so $m_1 \mid 4$. 

\noindent \textbf{Claim $(B)$: If $(i_3,j_3,m_3) = (4,5,3)$, then $F_{j_3} = K_3F_{i_1} = K_3\langle 1,v_1\rangle$ and $i_1 = 1$.} 
We have
\begin{gather}
\label{eqn:fj3}
\begin{aligned}
    F_{j_3} &\subseteq K_3F_{j_3} && \\
    &= K_3F_{i_1} && \text{because the inequality in \Cref{eqn:eqn-i1+1} is an equality} \\
\end{aligned}
\end{gather}

Now, we have
\begin{gather}
\label{eqn:eqn-i1+1}
\begin{aligned}
i_1+1 &= \dim_{L} F_{i_1} && \\
&= \frac{\dim_{L} K_3 F_{i_1}}{\dim_L K_3} && \text{because $F_{i_1} \subseteq K_1$ and $[K_1 : L] = 4$ and $[K_3 : L] = 3$ are coprime}\\ 
&\leq \frac{\dim_{L} K_3 F_{j_3}}{\dim_{L} K_3 } && \text{because $i_1 \leq j_3 = 5$}\\
&= \frac{2 \dim_{L} K_3}{\dim_{L} K_3 } && \text{by \Cref{lem:overflow}}\\
&= 2 &&
\end{aligned}
\end{gather}
Therefore, $i_1 = 1$, so the inequality in \Cref{eqn:eqn-i1+1} is an equality. 

Clearly, $\dim_{L}F_{j_3} = 6$ because $j_3 = 5$. As in \Cref{eqn:eqn-i1+1}, we have $\dim_{L}K_3F_{j_3} = 2 \dim_{L} K_3 = 6$, so the inequality in \Cref{eqn:fj3} is an equality. Therefore,
\begin{equation*}
    F_{j_3} = K_3F_{i_1} = K_3\langle 1,v_1\rangle
\end{equation*}
because $i_1 = 1$.

\noindent \textbf{Case 1: $\deg(v_1) = 4$.}
Suppose $\deg(v_1) = 4$. In this case $m = 4$ as well. 

\noindent \textbf{Claim $(1A)$: $(i_3,j_3,m_3) = (4,5,3)$.}
Suppose $i_3 < m_3$; then \Cref{lem:i-leq-m} implies that $F_{i_3} \subseteq K_{3}$; so $v_1 \in F_1 \subseteq F_{i_3} \subseteq K_{3}$, so $\deg(v_1) \mid m_{3}$. Now, looking at $\cL_3$ shows that $m_3 \in \{3,6\}$, which implies that $\deg(v_1) \mid 6$, which is a contradiction. 

Hence, we may suppose $i_3 > m_3$. By explicitly looking at $\cL_3$, we see that the only triple $(i_3,j_3,m_3)$ with  $i_3 > m_3$ is $(i_3,j_3,m_3) = (4,5,3)$. 

\noindent \textbf{Claim $(1B)$: $j_1 \in \{3,7\}$.}
Claim $(B)$ shows that $i_1 = 1$. By looking explicitly at $\cL_1$, we see that the only triples with $i_1 = 1$ and $m_1 = 4$ have $j_1 \in \{3,7\}$.

\noindent \textbf{Subcase $(a)$: $j_1 = 3$.} Suppose $j_1 = 3$. Then because $j_1 < m_1$, \Cref{lem:i-leq-m} implies that $F_{j_1} \subseteq K_1$. Because $\dim_{L}F_{j_1} = \dim_{L}K_1 = 4$, we have
\[
    F_{j_1} = K_1.
\]
Now, we have
\begin{gather}
\label{eqn:who-knows}
\begin{aligned}
    K_1 &= F_{j_1} && \\
    &\subseteq F_{j_3} && \text{because $j_1 \leq j_3$}\\
    &=K_3\langle 1,v_1\rangle && \text{by Claim $(B)$}.
\end{aligned}
\end{gather}
Because the latter is a $K_3$-vector space, \Cref{eqn:who-knows} implies that
\[
    K_3 K_1 \subseteq K_3\langle 1,v_1\rangle.
\]
However, $\dim_{L}K_3K_1 = \dim_{L}K_3 \dim_{L}K_1 = 3 \cdot 4 = 12$, and $\dim_{L}K_3\langle 1,v_1\rangle = 6$, which is a contradiction.

\noindent \textbf{Subcase $(b)$: $j_1 = 7$.} Suppose $j_1 = 7$. Then \Cref{lem:overflow} implies that
\[
    \dim_{L} K_1 F_{j_1} = 8.
\]
Because $F_{j_1} \subseteq K_1 F_{j_1}$ and $\dim_{L} F_{j_1} = 8$, we have $ K_1 F_{j_1} = F_{j_1}$, so $F_{j_1}$ is a $K_1$-vector space. Now, we also have
\begin{gather}
\begin{aligned}
K_3 &\subseteq F_{j_3} && \text{by Claim $(B)$} \\
&= F_{j_1} && \text{because $5 = j_3 \leq j_1 = 7$}. \\
\end{aligned}
\end{gather}
Now, because $F_{j_1}$ is a $K_1$-vector space containing $K_3$, it contains $K_1K_3$, which is a field of degree $12$. Thus we are done, as the dimension of $F_{j_1}$ is $8$.

\noindent \textbf{Case 2: $\deg(v_1) = 2$.}
Suppose $\deg(v_1) = 2$. Then
\[
F_1F_1 = L\langle 1,v\rangle L\langle 1,v\rangle = L\langle 1,v,v^2 \rangle = L\langle 1,v\rangle= F_1,
\]
and so $1 = T_{\cF}(1,1) = 1 < 2 = T_{(3,2,2)}(1,1)$. Moreover $(1,1)$ is a corner of $(3,2,2)$. Without loss of generality, we may suppose $i_1 = j_1 = 1$ and $m_1 = 2$, and so $K_1 = F_1$. Notice that $m_3 \in \{3,6\}$.

\noindent \textbf{Case 2a: $m_3 = 3$.} If $m_3 = 3$, then if $i_3 < m_3$ then \Cref{lem:i-leq-m} implies that $F_{i_3} \subseteq K_{3}$, and so $\deg(v_1) \mid K_{3} \mid 3$, which is a contradiction. Thus $i_3 > m_3$. Explicitly looking at $\cL_3$ shows that $(i_3,j_3,m_3) = (4,5,3)$. Now Claim $(B)$ implies that:
\begin{gather}
\begin{aligned}
F_{5} &= K_3F_{1} && \text{by Claim $(B)$} \\
&= K_3\langle 1,v\rangle && \text{by Claim (B)} \\
&= K_3 F_1 && \text{because $F_1 =  L\langle 1,v_1\rangle$} \\
&= K_3 K_1 && \text{because $K_1 = F_1$}.
\end{aligned}
\end{gather}
Therefore, $F_5$ is a number field of degree $6$. However, recall that by assumption $K_3 = \Stab(F_4F_5)$. But because $F_5$ is a field and $F_4 \subseteq F_5$, we have $F_4F_5 = F_5$, so $\Stab(F_4F_5) = F_5 \neq K_3$, which is a contradiction.

\noindent \textbf{Case 2b: $m_3 = 6$.}
Suppose $m_3 = 6$. For all $(i,j,m) \in \cL_2$, we have $i < m$, so \Cref{lem:i-leq-m} implies that $F_{i_2} \subseteq K_{2}$. Because $v_1 \in F_{i_2}$, we have $\deg(v_1) \mid \deg(K_2)$; looking explicitly at $\cL_2$ shows that $m_2 \in \{3,4\}$, so we must have $m_2 = 4$.

Looking explicitly at $\cL_2$ and $\cL_3$ shows that $i_2 < m_2$ and $i_3 < m_3$ and $2 \leq i_2,i_3$; applying \Cref{lem:i-leq-m} shows that $F_{i_2} \subseteq K_2$ and $F_{i_3} \subseteq K_3$. Because  $2 \leq i_2,i_3$, we have
\[
    F_2 \subseteq F_{i_2}\cap F_{i_3} \subseteq K_2 \cap K_3.
\]
Because $K_2$ is a degree $4$ field and $K_3$ is a degree $6$ field, their intersection has dimension at most $2$. However, $\dim_{ L}F_2 = 3$, so we have a contradiction.
\end{proof}

\section{Constructing orders with almost prescribed successive minima}
\label{sec:construction}

In \Cref{sec:joint-constraints} we proved joint constraints on the successive minima of orders in number fields arising from multiplication. In this section, we show that the constraints arising from multiplication are \emph{all} the constraints on successive minima by constructing orders with almost prescribed successive minima (\Cref{prop:direction-1-mult-spectrum}). We use this construction, along with the results of \Cref{sec:constraints} to provide a proof of \Cref{thm:sn-spectrum} and \Cref{thm:mult-spectrum}. 

\begin{proposition}
\label{prop:direction-1-mult-spectrum}
Let $K$ be a degree $n$ number field and let $\{1=v_0,\dots,v_{n-1}\}$ be a basis of $K$. Let $\cF$ be the flag given by $F_i = \QQ \la v_0,\dots,v_i\ra$ and let $\mathbf{x} \in P_{T_{\cF}}$ be a $\QQ$-point of the relative interior. Then there exists a family of orders $\{\cO_i\}_{i \in \ZZ_{\geq 1}} \subseteq K$ such that $\lim_{i \rightarrow \infty }\lvert \Disc(\cO_i)\rvert = \infty$ and $\lim_{i \rightarrow \infty} p_{\cO_i} = \mathbf{x}$.
\end{proposition}
\begin{proof}
Write the multiplication table of the $v_{i}$ as
\[
    v_iv_j = \sum_{k = 0}^{n-1}c_{ij}^k v_k.
\]
Set $x_0 \coloneqq 0$. Define $\cM$ to the set of $M \in \ZZ_{\geq 1}$ such that $M^{x_i + x_j - x_k}c_{ij}^k \in \ZZ$ for all $i,j,k \in [n]$.

\noindent \textbf{Claim: $\cM$ is infinite.}
Because $\mathbf{x}$ is a $\QQ$-point, $x_i + x_j - x_k$ is a rational number for all $i,j,k$. Therefore, to show that $\cM$ is an infinite set, it suffices to show that if $c_{ij}^k \neq 0$, then $x_i + x_j - x_k \geq 0$ (equivalently, $x_k \leq x_i + x_j$). 

Now, if $c_{ij}^k \neq 0$, then $F_i F_j \not \subseteq F_{k-1}$, so $T_{\cF}(i,j) \geq k$. Therefore, $P_{T_{\cF}}$ is contained in the linear half-space given by $x_{T_{\cF}(i,j)} \leq x_i + x_j$. Because $k \leq T_{\cF}(i,j)$, $P_{T_{\cF}}$ is contained in the linear half-space given by $x_k \leq x_i + x_j$, so the claim is proven. 

For $M \in \cM$, define the free $\ZZ$-module:
\[
	\cO_M \coloneqq \ZZ \la 1 = M^{x_0}v_0, M^{x_1}v_1,\dots, M^{x_{n-1}} v_{n-1} \ra.
\]

\noindent \textbf{Claim: $\cO_M$ is a ring.} 
We have:
\[
	(M^{x_i}v_i)(M^{x_j}v_j) = \sum_{k \in [n]} M^{x_i + x_j - x_k}c_{ij}^k(v_iv_j) (M^{x_k}v_k).
\]
Now, by assumption, $M^{x_i + x_j - x_k}c_{ij}^k \in \ZZ$.

\noindent \textbf{Claim: $\lim_{M \rightarrow \infty}p_{\cO_M} = \mathbf{x}$.}
We have:
\[
    \Disc(\cO_M) = \Disc(\ZZ\langle v_0,\dots,v_{n-1}\rangle)M^{2(\sum_{i}x_i)} = \Disc(\ZZ\langle v_0,\dots,v_{n-1}\rangle)M
\]
Thus, we obtain:
\begin{align*}
M^{1/2} &\asymp_{v_1,\dots,v_{n-1}} \lvert \Disc(\cO_M) \rvert^{1/2} &&  \\
&\asymp_n \prod_{i = 1}^{n-1} \lambda_i(\cO_M) && \text{by Minkowski's second theorem} \\
&\ll_{v_1,\dots,v_{n-1}} \prod_{i = 1}^{n-1} M^{x_i} && \text{because } \lambda_i(\cO_M) \ll_{v_1,\dots,v_{n-1}} M^{x_i} \\
&= M^{1/2} && \text{because } x_1 + \dots + x_{n-1} = 1/2
\end{align*}
This implies that:
\[
	\lambda_i(\cO_M) \asymp_{v_1,\dots,v_{n-1}} M^{x_i}.
\]
Therefore,
\[
	\lim_{M \rightarrow \infty}\log_{\lvert \Disc(\cO_M) \rvert}\lambda_i(\cO_M) = \lim_{M \in \cM}\log_{M}M^{x_i} = x_i.
\]
\end{proof}

\begin{proof}[Proof of \Cref{thm:mult-spectrum}]
Let $\cO$ be an order in a degree $n$ number field with Galois group $G$. Because $1 \leq \lambda_1 \leq \dots \leq \lambda_{n-1}$ and $\prod_i \lambda_i \asymp_n \lvert \Disc(\cO) \rvert^{1/2}$, we have that
\begin{equation}
\label{eqn:basic-constraints}
    \Spe(\Sigma(G)) \subseteq \{\mathbf{x} \in \RR^{n-1} : \sum_{i=1}^{n-1} x_i = 1/2 \text{ and } 0 \leq x_1 \leq \dots \leq x_{n-1}\}.
\end{equation}

Now let $\{v_0=1, v_1, v_2,\dots,v_{n-1}\}$ be a Minkowski reduced basis of $\cO$, let $\cF$ be the corresponding flag, and let $T_{\cF}$ be the corresponding flag type. \Cref{prop:flag-bound} shows that $\lambda_{T_{\cF}(i,j)} \ll_n \lambda_i \lambda_j$ for all $1 \leq i,j < n$. Therefore, 
\begin{equation}
\label{eqn:nontriv-constraints}
\Spe(\Sigma(G)) \subseteq \{\mathbf{x} \in \RR^{n-1} : x_{T_{\cF}(i,j) }\leq x_i + x_j\}.
\end{equation}
as $\cF$ ranges across flags of degree $n$ number fields with Galois group $G$. Combining \Cref{eqn:basic-constraints} and \Cref{eqn:nontriv-constraints}, we see that 
\begin{equation}
\label{eqn:one-inclusion}
\Spe(\Sigma(G)) \subseteq \bigcup_{\cF} P_{T_{\cF}}.
\end{equation}

Conversely, let $\cF$ be a flag of a degree $n$ extension $K$. Choose a basis $\{v_0=1, v_1, v_2,\dots,v_{n-1}\}$ of $K$ such that $F_i = \QQ\la v_0,\dots,v_i\ra$. Then \Cref{prop:direction-1-mult-spectrum} shows that
\[
    \QQ^{n-1} \cap P_{T_{\cF}} \subseteq \Spe(\Sigma(G)).
\]

Now, $\Spe(\Sigma(G))$ is defined to be the set of limit points of a multiset; hence, it is closed. Therefore, 
\[
    \overline{\QQ^{n-1} \cap P_{T_{\cF}}} = P_{T_{\cF}} \subseteq \Spe(\Sigma(G)).
\]
As we range across all flags $\cF$ of degree $n$ extensions with Galois group $G$, we obtain 
\begin{equation}
\label{eqn:second-inclusion}
\bigcup_{\cF} P_{T_{\cF}} \subseteq \Spe(\Sigma(G)).
\end{equation}
Combining \Cref{eqn:one-inclusion} and \Cref{eqn:second-inclusion}, we get
\[
\bigcup_{\cF} P_{T_{\cF}} = \Spe(\Sigma(G)).
\]

\end{proof}

\subsection{Computing $\Spe(\Sigma(S_n))$}
Using \Cref{thm:mult-spectrum}, we now compute $\Spe(\Sigma(S_n))$. We'll need the following lemma, which shows that the polytope $\Len_{(n)}$ is equal to $P_{T_{\cF}}$ for some flag $\cF$. 

\begin{lemma}
\label{lem:len-n}
Let $K$ be any degree $n$ number field. Choose $\alpha \in K$ such that $\QQ(\alpha) = K$. Let $\cF$ be the flag such that $F_i = \QQ \la 1,\alpha,\dots,\alpha^i \ra$. Then $P_{T_{\cF}} = \Len_{(n)}$.
\end{lemma}
\begin{proof}
We see that
\[
    F_iF_j = \QQ \la 1,\alpha,\dots,\alpha^i \ra\QQ \la 1,\alpha,\dots,\alpha^j \ra = \QQ \la 1,\alpha,\dots,\alpha^{i+j} \ra = \QQ \la 1,\alpha,\dots,\alpha^{\min(n-1,i+j)} \ra = F_{\min( n-1,i+j)}
\]
Therefore, $T_{\cF}(i,j) = \min(n-1, i+j)$ for all $i,j$. Therefore, $P_{T_{\cF}}$ is defined by the inequalities:
\begin{itemize}
    \item $\sum_{i = 1}^{n-1}x_i = 1/2$;
    \item $0 \leq x_1 \leq \dots \leq x_{n-1}$;
    \item and $x_{\min(n-1,i+j)} \leq x_i + x_j$ for all $1 \leq i,j < n$.
\end{itemize}
By removing extraneous inequalities, we see that $P_{T_{\cF}}$ is defined by the inequalities:
\begin{itemize}
    \item $\sum_{i = 1}^{n-1}x_i = 1/2$;
    \item $0 \leq x_1 \leq \dots \leq x_{n-1}$;
    \item and $x_{i+j} \leq x_i + x_j$ for all $1 \leq i,j < i +j < n$.
\end{itemize}
These are precisely the inequalities defining $\Len_{(n)}$.
\end{proof}

\begin{proof}[Proof of \Cref{thm:sn-spectrum}] 
\Cref{thm:bl} implies that $\Spe(\Sigma(S_n))$ is contained in $\{\mathbf{x} \in \RR^{n-1} : x_{i+j} \leq x_i + x_j \; \forall i,j\}$. Moreover, because $1 \leq \lambda_1 \leq \dots \leq \lambda_{n-1}$ and $\prod_i \lambda_i \asymp_n \lvert \Disc(\cO) \rvert^{1/2}$, we have that
\begin{equation*}
    \Spe(\Sigma(G)) \subseteq \{\mathbf{x} \in \RR^{n-1} : \sum_{i=1}^{n-1} x_i = 1/2 \text{ and } 0 \leq x_1 \leq \dots \leq x_{n-1}\}.
\end{equation*}
Together, these two containments imply that $\Spe(\Sigma(S_n)) \subseteq \Len_{(n)}$.

Conversely, let $K$ be any degree $n$ number field with Galois group $S_n$. Choose $\alpha \in K$ such that $\QQ(\alpha) = K$. Let $\cF$ be the flag such that $F_i = \QQ \la 1,\alpha,\dots,\alpha^i \ra$ for all $i$. \Cref{thm:mult-spectrum} shows that $P_{T_{\cF}} \subseteq \Spe(\Sigma(S_n))$, and \Cref{lem:len-n} shows that $P_{T_{\cF}} = \Len_{(n)}$, so $\Len_{(n)} \subseteq \Spe(\Sigma(S_n))$. 
\end{proof}

\section{Computing $\Spe(\Sigma(S_n))$ when $n$ is a prime power, a product of $2$ primes, or $12$}
\label{sec:lenstra-spectrum-good}

In \Cref{sec:joint-constraints}, we explicitly described the flag types which occur from flags when $n$ is a prime power, a product of $2$ primes, or $12$. In \Cref{sec:construction} we explicitly described the successive minima spectrum in terms of flag types. In this (short) section, we combine these two results to more explicitly describe the successive minima spectrum when $n$ is a prime power, a product of $2$ primes, or $12$.

Namely, we prove \Cref{conj:lenstra} when $n$ is a prime power, a product of $2$ primes, or $12$. In particular, we prove \Cref{thm:lenstra-spectrum} in these cases.

\begin{proof}[Proof of \Cref{thm:lenstra-spectrum} when $n$ is a prime power, a product of $2$ primes, or $12$]
Let  $n$ be a prime power, a product of $2$ primes, or $12$. Then \Cref{prop:good-n-direction-1} and \Cref{prop:good-n-direction-2} together imply that $\Spe(\Sigma_n) = \cup_{\mfT} \Len_{\mfT}$.
\end{proof}

\begin{proposition}
\label{prop:good-n-direction-1}
Suppose $n$ is a prime power, $12$, or a product of two primes. Then 
\[
    \Spe(\Sigma_n) \subseteq \bigcup_{\mfT} \Len_{\mfT}.
\]
\end{proposition}
\begin{proof}
\Cref{thm:containment} along with the explicit description of the Lenstra polytopes given in \Cref{def:lenstra-polytope}, proves the proposition.
\end{proof}

To show $\bigcup_{\mfT} \Len_{\mfT} \subseteq \Spe(\Sigma_n)$, we'll need the following crucial lemma. We delay the proof to \Cref{subsec:len-T}. 

\begin{lemma}
\label{lem:len-T}
Let $\mfT = (n_1,\dots,n_t)$ be a tower type. Let $\alpha_1,\dots,\alpha_{t} \in \overline{\QQ}$ be elements such that $\deg(\alpha_i) = n_i$, the field $\QQ(\alpha_i)$ has no nontrivial proper subfields, and the compositum $\QQ(\alpha_1,\dots,\alpha_t)$ has degree $n$. Set $K = \QQ(\alpha_1,\dots,\alpha_t)$. For $1 \leq j < n$, write $j$ in mixed radix notation with respect to $\mfT$ as
\[
    j = j_1 + j_2 n_1 + i_3 (n_1n_2) + \dots + j_t(n_1\dots n_{t-1}).
\]
Define a basis $\{1 = v_0,\dots,v_{n-1}\}$ of $K$ by setting $v_j \coloneqq \prod_{\ell = 1}^t \alpha_\ell^{j_\ell}$. Let $\cF$ be the corresponding flag. Then, $P_{T_{\cF}} = \Len_{\mfT}$.
\end{lemma}

\begin{proposition}
\label{prop:good-n-direction-2}
For all $n$, we have
\[
    \bigcup_{\mfT} \Len_{\mfT} \subseteq \Spe(\Sigma_n).
\]
\end{proposition}
\begin{proof}
\Cref{thm:mult-spectrum} states that
\[
    \Spe(\Sigma_n) = \bigcup_{\cF} P_{T_{\cF}}
\]
as $\cF$ ranges across flags in degree $n$ number fields. Now, \Cref{lem:len-T} shows that for every tower type $\mfT$, there is a flag $\cF$ such that $P_{T_{\cF}} = \Len_{\mfT}$. Therefore, 
\[
    \bigcup_{\mfT} \Len_{\mfT} \subseteq \Spe(\Sigma_n).
\]
\end{proof}

\subsection{Showing that for every $\mfT$, there exists a flag $\cF$ so that $\Len_{\mfT} = P_{T_{\cF}}$}
\label{subsec:len-T}

\begin{proof}[Proof of \Cref{lem:len-T}]
From the definition of $P_{T_{\cF}}$, we see that $P_{T_{\cF}}$ is defined by the inequalities:
\begin{itemize}
    \item $\sum_{i = 1}^{n-1}x_i = 1/2$;
    \item $0 \leq x_1 \leq \dots \leq x_{n-1}$;
    \item and $x_{T_{\cF}(i,j)} \leq x_i + x_j$ for all $1 \leq i,j < n$. 
\end{itemize}
We now explicitly describe the third inequality. Choose $1 \leq i,j < n$. Write $i,j$ in mixed radix notation as above. Set 
\[
     k \coloneqq \min(n_1 - 1, i_1 + j_1) + \min(n_2 - 1, i_2 + j_2) n_1 +  \dots + \min(n_t - 1, i_t + j_t)(n_1\dots n_{t-1}).
\]

\noindent \textbf{Claim: $v_iv_j \in F_{k} \setminus F_{k-1}$.}
We have 
\[
    v_i v_j = \prod_{\ell = 1}^t \alpha_\ell^{i_\ell} \prod_{\ell = 1}^t \alpha_\ell^{j_\ell} = \prod_{\ell = 1}^t \alpha_\ell^{i_\ell + j_\ell}. 
\]
If $i_\ell + j_\ell < n_\ell$ for all $1 \leq \ell \leq t$, then $v_iv_j = v_{i+j}$ and $i + j = k$. Therefore $v_i v_j \in F_{k} \setminus F_{k - 1}$.


On the other hand, let $S = \{\ell : i_\ell + j_\ell \geq n_\ell\}$. Then we may write:
\begin{align*}
 v_i v_j &= \prod_{\ell \notin S} \alpha_\ell^{i_\ell + j_\ell} \prod_{\ell \in S} \alpha_\ell^{i_\ell + j_\ell} \\
 &\in \prod_{\ell \notin S} \alpha_\ell^{i_\ell + j_\ell} \prod_{\ell \in S} \QQ \la 1,
\alpha_k,\dots,\alpha_{\ell}^{n_\ell - 1}\ra \\
&\subseteq F_{k}
\end{align*}

So, $v_iv_j \in F_{K}$. Because $\QQ(\alpha_i)$ has no nontrivial proper subfields, the coefficient of $\alpha_\ell^{n_i - 1}$ is nonzero in the expansion of $\alpha_\ell^{i_\ell + j_\ell}$ for all $\ell \in S$. Therefore, $v_i v_j \notin F_{k - 1}$.

\noindent \textbf{Describing $P_{T_{\cF}}$.} Letting $L$ be as above, we see that $P_{T_{\cF}}$ is defined by the inequalities:
\begin{itemize}
    \item $\sum_{i = 1}^{n-1}x_i = 1/2$;
    \item $0 \leq x_1 \leq \dots \leq x_{n-1}$;
    \item and $x_{k} \leq x_i + x_j$ for all $1 \leq i,j < n$. 
\end{itemize}
Moreover, observe that $\ell = i + j$ if $i + j$ does not overflow modulo $\mfT$, and $\ell < i+j$ if $\ell$ overflows modulo $\mfT$. Removing extraneous inequalities, we see that $P_{T_{\cF}}$ is defined by the inequalities:
\begin{itemize}
    \item $\sum_{i = 1}^{n-1}x_i = 1/2$;
    \item $0 \leq x_1 \leq \dots \leq x_{n-1}$;
    \item and $x_{i+j} \leq x_i + x_j$ if $i+j$ does not overflow modulo $\mfT$. 
\end{itemize}
Now, these are precisely the inequalities defining $\Len_{\mfT}$. 
\end{proof}

\section{Proving $\Spe(\Sigma(S_n)) \neq \cup_{\mfT}\Len_{\mfT}$ when $n$ is not a prime power, a product of $2$ primes, or $12$}
\label{sec:lenstra-spectrum-bad}

In this section, we give a proof of \Cref{thm:lenstra-spectrum} in the case when $n$ is not a prime power, a product of $2$ primes, or $12$. Combined with the results of \Cref{sec:lenstra-spectrum-good}, this completes the proof of \Cref{thm:lenstra-spectrum}.

\begin{proposition}
\label{prop:counterexamples}
Suppose $n$ is not a prime power, $12$, or a product of two primes. Then:
\[
    \Spe(\Sigma_n) \not \subseteq \bigcup_{\mfT} \Len_{\mfT}.
\]
\end{proposition}
\begin{proof}
\Cref{thm:mult-spectrum} says that
\[
    \Spe(\Sigma_n) = \bigcup_{\cF} P_{T_{\cF}}
\]
as $\cF$ ranges over flags in degree $n$ fields. Therefore, to prove the proposition, it suffices to show that there exists a flag $\cF$ such that 
\begin{equation}
\label{eqn:big-ptf}
P_{T_{\cF}} \not \subseteq \bigcup_{\mfT} \Len_{\mfT}.
\end{equation}

By \Cref{lem:not-complete}, the existence of such a flag for degree $m$ implies the existence of such a flag for degree $n$, where here $m \mid n$. Therefore, it suffices to show the existence of such a flag when:
\begin{enumerate}
\item $n = p^2q$ for two distinct odd primes $p$ and $q$ with $p < q$, in which case \cref{lem:p2q} provides a proof;
\item $n = pqr$ for three primes $p$, $q$, and $r$ with $p < q \leq r$, in which case \cref{lem:pqr} provides a proof;
\item $n = 4p$ for a prime $p \neq 2,3$, in which case \cref{lem:4p} provides a proof;
\item or $n = 24$, in which case \cref{lem:24} provides a proof.
\end{enumerate}
\end{proof}

\begin{definition}
\label{def:cone}
Given a set $S \subseteq \RR^k$,  we say \emph{the cone over $S$} is 
\[
    \Cone(S) \coloneqq \{\alpha \mathbf{x} : \alpha \in \RR_{\geq 0}, \; \mathbf{x} \in S\}.
\]
\end{definition}

\begin{proposition}
\label{prop:dim}
Let $T$ be any flag type. Then the set $P_T$ is a bounded polytope of dimension $n-2$.
\end{proposition}
\begin{proof}
Note that $P_T$ lies in the hyperplane in $\RR^{n-1}$ whose coordinates sum to $1/2$. Thus to showing that $P_T$ has dimension $n-2$ is equivalent to showing that the cone over $P_T$ contains $n-1$ linearly independent vectors. For $1 \leq \ell \leq n-1$, define $\mathbf{x}^\ell = (x_1^\ell,\dots,x_{n-1}^\ell) \in \RR^{n-1}$ by
\begin{align*}
x_1^\ell = \dots = x_\ell^\ell &= 1 \\
x_{\ell+1}^\ell = \dots = x_{n-1}^\ell &= 2.
\end{align*}
Clearly, $0 \leq x_1^{\ell} \leq \dots \leq x_{n-1}^{\ell}$ and for all $1 \leq i,j,k < n$, we have $x_k^\ell \leq x_i^\ell + x_j^\ell$, so $\mathbf{x}^\ell$ is contained in the cone over $P_T$. Consider the matrix whose columns are the $\mathbf{x}^{\ell}$:
\[
\begin{bmatrix}
    x_1^{1}       & x_1^2 & x_{1}^3 & \dots & x_{1}^{n-1} \\
    x_{2}^2       & x_{2}^2 & x_{2}^3 & \dots & x_{2}^{n-1} \\
    \hdotsfor{5} \\
    x_{n-1}^1       & x_{n-1}^2 & x_{n-1}^3 & \dots & x_{n-1}^{n-1}.
\end{bmatrix}
\]
Modulo $2$, the matrix is equal to the upper triangular matrix
\[
\begin{bmatrix}
    1       & 1 & 1& \dots & 1 \\
    0       & 1 & 1 & \dots & 1 \\
    \hdotsfor{5} \\
    0       & 0 & 0 & \dots & 1.
\end{bmatrix} 
\]
which visibly has nonzero determinant. Thus the $\mathbf{x}^\ell$ form a set of $n-1$ linearly independent vectors in the cone over $P_T$. 
\end{proof}

\begin{lemma}
\label{lem:not-complete}
Let $n,m \in \ZZ_{>1}$ be integers such that $m \mid n$.
If there exists a flag $\cF$ of a degree $m$ number field such that 
\[
	P_{T_{\cF}} \not\subseteq \bigcup_{\mfT}\Len_{\mfT}
\]
where $\mfT$ ranges across tower types of degree $m$, then there exists a flag $\cF'$ of a degree $n$ number field such that 
\[
	P_{T_{\cF'}} \not\subseteq \bigcup_{\mfT'}\Len_{\mfT'}
\]
where $\mfT'$ ranges across tower types of degree $n$.
\end{lemma}
\begin{proof}
By induction, it suffices to assume $n = pm$ for $p$ a prime. Let $K$ denote the degree $m$ number field containing the flag $\cF$ and let $1 = v_0,\dots,v_{m-1} \in K$ be such that $F_i = \QQ\la v_0,\dots,v_{i} \ra$. Let $L$ be a degree $p$ extension of $K$, and let $\alpha \in L$ be such that $L=K(\alpha)$. Define the sequence $\{1 = v'_0,\dots,v'_{n-1} \}$  by $v'_i = v_{i_2}\alpha^{i_1}$ for $i = i_1 + i_2 p$ in mixed radix notation with respect to $(p,m)$. Define a flag $\cF' = \{F'_i\}_{i \in [n]}$ of $L$ by setting $F'_i \coloneqq \QQ \la v'_0,\dots,v'_i \ra$.

By \Cref{prop:dim}, the polytope $P_{T_{\cF}}$ has dimension $m-2$. Because $P_{T_{\cF}} \setminus \bigcup_{\mfT}\Len_{\mfT}$ is nonempty by assumption, it must also have dimension $m - 2$. Therefore, $P_{T_{\cF}} \setminus \bigcup_{\mfT}\Len_{\mfT}$ is full-dimensional inside the hyperplane $\{\mathbf{x} \in \RR^{m-1} : \sum_{i = 1}^{m-1}x_i = 1/2\}$. As a result, the set
\[
	S \coloneqq \big(P_{T_{\cF}} \setminus \bigcup_{\mfT}\Len_{\mfT}\big) \cap \{\mathbf{x} \in \RR^{m-1} : x_i \neq x_{i+1} \; \forall  \; 1 \leq i < m-1 \}
\] 
is nonempty. 

Choose $\mathbf{x} = (x_1,\dots,x_{m-1}) \in S$ and set $\epsilon \coloneqq \min_{1 \leq i < m-1}\{x_{i+1} - x_i\}$. By definition we have $\epsilon \neq 0$. Define the point $\mathbf{x}' = (x_1',\dots,x_{n-1}') \in \RR^{n-1}$ as follows. For every $1 \leq i < n$, write $i = i_1 + i_2p$ in mixed radix notation with respect to $(p,m)$ and set
\[
	x_i' \coloneqq \epsilon \frac{i_1}{2p} + x_{i_2}.
\]

\noindent \textbf{Claim: $\mathbf{x}' \in \Cone(P_{T_{\cF'}})$.} It suffices to show that $0 \leq x_1' \leq \dots \leq x'_{n-1}$ and that $x'_{T_{\cF'}(i,j)} \leq x'_i + x'_j$ for all $1 \leq i,j < n$. 

For $1 \leq i < n-1$, write $i = i_1 + i_2 p$ in mixed radix notation with respect to $(p,m)$. If $i_1 \neq p-1$, then
\begin{equation}
\label{eqn:induction-1}
x_{i+1}' = \epsilon \frac{i_1 + 1}{2p} + x_{i_2} \geq  \epsilon \frac{i_1}{2p} + x_{i_2} = x_i'.
\end{equation}
If $i_1 = p-1$ then 
\begin{equation}
\label{eqn:induction-2}
x_{i+1}' = x_{i_2+1} \geq  \epsilon + x_{i_2} \geq \epsilon \frac{i_1}{2p} + x_{i_2} = x_i'.
\end{equation}
Combining \Cref{eqn:induction-1} and \Cref{eqn:induction-2}, we see that  $0 \leq x_1' \leq \dots \leq x'_{n-1}$.

For any $1 \leq i,j < n$, we now show that $x'_{T_{\cF'}(i,j)} \leq x'_i + x'_j$. Let $k = T_{\cF'}(i,j)$. Write
\begin{align*}
i &= i_1 + i_2 p \\
j &= j_1 + j_2 p \\
k &= k_1 + k_2 p
\end{align*}
in mixed radix notation with respect to $(p,m)$. The explicit description of the $v'_\ell$ shows that $k_1 \leq i_1+j_1$ and $T_{\cF}(i_2,j_2) = k_2$, so $x_{k_2}\leq x_{i_2}+x_{j_2}$. Then
\begin{align*}
x'_k &= \epsilon \frac{k_1}{2p} + x_{k_2} \\
&\leq \epsilon \frac{k_1}{2p} + x_{k_2}\\
&\leq \epsilon\bigg(\frac{i_1}{2p} + \frac{j_1}{2p}\bigg) + x_{i_2} + x_{j_2} \\
&\leq \epsilon \frac{i_1}{2p} + x_{i_2} + \epsilon \frac{j_1}{2p} + x_{j_2} \\
&\leq x'_i + x'_j.
\end{align*}

\noindent \textbf{Claim: $\mathbf{x}' \notin \cup_{\mfT'}\Cone(\Len_{\mfT'})$ where $\mfT'$ ranges across tower types of degree $n$.}
First, notice that
\[
    \bigcup_{\mfT'}\Len_{\mfT'} = \bigcup_{(p_1,\dots,p_t)}\Len_{(p_1,\dots,p_t)}
\]
where $(p_1,\dots,p_t)$ ranges across all tuples with prime entries such that $\prod_i p_i = n$. Fix such a tuple $(p_1,\dots,p_t)$.

If $p_1 \neq p$ then
\begin{equation}
\label{eqn:induction-3}
x_1'+x_{p-1}' = \frac{\epsilon}{2} < x_1 = x'_p.
\end{equation}
Now, by definition,
\[
    \Len_{(p_1,\dots,p_t)} \subseteq \{\mathbf{y} \in \RR^{n-1} : y_p \leq y_1 + y_{p-1}\}
\]
because $1 + (p-1)$ does not overflow modulo $(p_1,\dots,p_t)$. Therefore \Cref{eqn:induction-3} implies that
\[
    \mathbf{x}' \notin \Cone(\Len_{(p_1,\dots,p_t)}).
\]

If $p_1 = p$ then $p_2\dots p_t = m$. Because $\mathbf{x} \notin \Cone(\Len_{(p_2,\dots,p_t)})$ by assumption, we can choose $1 \leq i \leq j < i+j < m$ such that $i+j$ does not overflow modulo $(p_2,\dots,p_t)$ and $x_{i+j} > x_i + x_j$. Note that:
\begin{equation}
\label{eqn:induction-4}
x_{pi+pj}' = x_{i+j} > x_{i} + x_{j} = x_{pi}' + x_{pj}'.
\end{equation}
Now, by definition,
\[
    \Len_{(p_1,\dots,p_t)} \subseteq \{\mathbf{y} \in \RR^{n-1} : y_{pi + pj} \leq y_{pi} + y_{pj}\}
\]
because $pi + pj$ does not overflow modulo $(p_1,\dots,p_t)$. Therefore, \Cref{eqn:induction-4} implies that
\[
    \mathbf{x}' \notin \Cone(\Len_{(p_1,\dots,p_t)}).
\]

\noindent \textbf{Completing the proof.} Let $\mathbf{x}''$ be the point obtained by scaling $\mathbf{x}$ so the coordinates sum to $1/2$. Both claims together imply that $\mathbf{x}'' \in P_{T_{\cF'}} \setminus \cup_{\mfT'}\Len_{\mfT'}$.
\end{proof}

\begin{proposition}
\label{lem:p2q}
Let $p$ and $q$ be two distinct odd primes with $p < q$ and let $n = p^2q$. Then there exists a flag $\cF$ of a degree $n$ number field such that 
\[
	P_{T_{\cF}} \not \subseteq \Len_{(p,p,q)} \cup \Len_{(p,q,p)} \cup \Len_{(q,p,p)}.
\]
\end{proposition}
\begin{proof}

It suffices to show that there exists a flag $\cF$ and a point
\[
	\mathbf{x} \in \Cone(P_{T_{\cF}}) \setminus \Cone( \Len_{(p,p,q)} \cup \Len_{(p,q,p)} \cup \Len_{(q,p,p)}).
\] 

\noindent \textbf{Part A: defining the flag $\cF$.} Choose $e_1,e_2,e_3 \in \overline{\QQ}$ such that:
\begin{itemize}
    \item $e_1$ and $e_2$ have degree $p$;
    \item $e_3$ has degree $q$;
    \item and the compositum $K = \QQ(e_1,e_2,e_3)$ has degree $p^2q$.
\end{itemize}
Define a basis $\{1=v_0,\dots,v_{n - 1}\}$ of $K$ as follows. For $0 \leq i < q$, set $v_i \coloneqq e_1^{i}$. For $q \leq i < n$ and $1 \leq i' < n$, write $i' = i'_1 + i'_2p + i_3'pq$ in mixed radix notation with respect to $(p,q,p)$. 
Inductively define $v_i$ as follows. Choose $i'$ minimal such that $e_2^{i'_1}e_1^{i'_2}e_3^{i'_3} \notin \{v_0,\dots,v_{i-1}\}$, and set $v_i \coloneqq e_2^{i'_1}e_1^{i'_2}e_3^{i'_3}$. Define a flag $\cF = \{F_i\}_{i \in [n]}$ by $F_i \coloneqq \QQ \la v_0,\dots,v_i \ra$.

\noindent \textbf{Part B: explicit description of $\cF$.}
We first make some explicit descriptions of the flag $\cF$. Recall that for an element $\alpha \in K$, $\QQ(\alpha)$ refers to the field generated by $\alpha$. For a field $L \subseteq K$, $L \la \alpha\ra$ refers to the $L$-vector space generated by $\alpha$. For two $L$-vector spaces $A,B \subseteq K$, the sum $A + B = \{a + b : a \in A, b \in B\}$.

First, it follows immediately from the definitions that
\begin{equation}
\label{eqn:p2q-f1}
F_1  = \QQ \la  1, e_1 \ra
\end{equation}
\begin{equation}
\label{eqn:p2q-fq-1}
F_{q-1} = \QQ(e_1).
\end{equation}

\noindent \textbf{Claim: $\{v_0,\dots,v_{pq-1}\} = \{e_2^{i'_1}e_1^{i'_2}e_3^{i'_3} : i' < pq\}$.} It is clear from definition that $\{v_{q},\dots,v_{pq-1}\} \subseteq \{e_2^{i'_1}e_1^{i'_2}e_3^{i'_3} : i' < pq\}$. For every $0 \leq j < q$,
\[
    v_j = e_1^j = e_2^{0}e_1^{j}e_3^{0} = e_2^{i'_1}e_1^{i'_2}e_3^{i'_3}
\]
where $i' = 0 + jp + 0(pq)$. Because $j < q$, we have $i' < pq$.

Next, the claim above shows that 
\begin{equation}
\label{eqn:p2q-fpq}
F_{pq} = \QQ(e_1,e_2).
\end{equation}
Explicit computation shows that $v_{pq} = e_3$ and $v_{pq + 1} = e_2e_3$. Therefore:
\begin{equation}
\label{eqn:p2q-fpq+1}
F_{pq+1} = \QQ(e_1,e_2) + \QQ \la  e_3, e_2e_3 \ra
\end{equation}

The claim above implies that for $q \leq i < n$, we have $v_i = e_2^{i_1}e_1^{i_2}e_3^{i_3}$ when $i$ is in mixed radix notation with respect to $(p,q,p)$. Therefore:
\begin{equation}
\label{eqn:p2q-fpq+p-1}
F_{pq+p-1} = \QQ(e_1,e_2) + \QQ(e_2) \QQ \la  e_3\ra = \QQ(e_2)( \QQ(e_1) + \QQ \la e_3 \ra).
\end{equation}

Similarly, 
\begin{equation}
\label{eqn:p2q-2pq+p-1}
F_{2pq+p-1} = \QQ(e_1, e_2) \QQ \la 1, e_3 \ra + \QQ(e_2) \QQ \la e_3^2 \ra
\end{equation}
It is easy to see that:
\begin{equation}
\label{eqn:f1fq-1}
F_1 F_{q-1} = F_{q-1} 
\end{equation}
and
\begin{gather}
\begin{aligned}
\label{eqn:fpq+1-pq+p-1}
F_{pq+1}F_{pq+p-1} &= (\QQ(e_1, e_2) + \QQ\la e_3, e_2e_3 \ra)(\QQ(e_2) (\QQ( e_1 ) + \QQ\la e_3 \ra))\\
&= \QQ( e_1, e_2) \QQ \la 1, e_3 \ra + \QQ(e_2)\QQ \la e_3^2 \ra\ \\
&= F_{2pq+p-1}.
\end{aligned}
\end{gather}

\noindent \textbf{Part C: showing there exists $\mathbf{x} \in \Cone(P_{T_{\cF}})$ such that $x_q > x_1 + x_{q-1}$ and $x_{2pq+p} > x_{pq + 1} + x_{pq + p - 1}$.}
For $1 \leq i < q$ set $x_i \coloneqq \frac{i}{2q}$. For $q \leq i < pq$ set $x_i = 1$. For $pq \leq i < p^2q$ write $i = i_1 + i_2p + i_3pq$ in mixed radix notation with respect to $(p,q,p)$ and set $x_i = i_1 \frac{1}{4pq} + i_2 \frac{1}{2q} + i_3$. 

It is easy to see that:
\[
	x_q = 1 > \frac{1}{2q} + \frac{q-1}{2q} = x_1 + x_{q-1}
\]
and
\[
	x_{2pq+p} = \frac{1}{2q} + 2 > \bigg (\frac{1}{4pq} + 1 \bigg) + \bigg (\frac{p-1}{4pq} + 1 \bigg) = x_{pq+1} + x_{pq + p - 1}.
\]
So it remains to show that $\mathbf{x} \in \Cone(P_{T_{\cF}})$.

\noindent \textbf{Part C.1: showing that $0 \leq x_1 \leq \dots x_{n-1}$.}
If $1\leq i < q-1$, then
\[
	x_{i+1} = \frac{i+1}{2q} \geq \frac{i}{2q} = x_i.
\]
Note also that
\[
	x_q = 1 > \frac{q-1}{2q} = x_{q-1}
\]
If $q\leq i < pq-1$, 
\[
	x_{i+1} = 1 = x_i.
\]
If $pq\leq i < n-1$ then write $i+1 = (i+1)_1 + (i+1)_2p + (i+1)_3pq$ in mixed radix notation with respect to $(p,q,p)$ as well. Then if $i_1 = p-1$ and $i_2 = q-1$ then $(i+1)_1 = 0$ and $(i+1)_2 = 0$ and $(i+1)_3 = i_3 + 1$. Then
\begin{align*}
x_{i+1} &= (i+1)_1 \frac{1}{4pq} + (i+1)_2 \frac{1}{2q} + (i+1)_3 \\
&= (i_3 + 1) \\
& > i_1 \frac{1}{4pq} + i_2 \frac{1}{2q} + i_3 \\
&= x_i.
\end{align*}
Instead if $i_1 = p-1$ and $i_2 \neq q-1$ then $(i+1)_1 = 0$ and $(i+1)_2 = i_2 + 1$ and $(i+1)_3 = i_3$. Then
\begin{align*}
x_{i+1} &= (i+1)_1 \frac{1}{4pq} + (i+1)_2 \frac{1}{2q} + (i+1)_3 \\
&= (i_2+1) \frac{1}{2q} + i_3 \\
&> i_1 \frac{1}{4pq} + i_2 \frac{1}{2q} + i_3 \\
&= x_i.
\end{align*}
Finally, if $i_1 \neq p-1$ then $(i+1)_1 = i_1+1$ and $(i+1)_2 = i_2$ and $(i+1)_3 = i_3$. Then
\begin{align*}
x_{i+1} &= (i+1)_1 \frac{1}{4pq} + (i+1)_2 \frac{1}{2q} + (i+1)_3 \\
&= (i_1+1) \frac{1}{4pq} + i_2\frac{1}{2q} + i_3 \\
&> i_1 \frac{1}{4pq} + i_2 \frac{1}{2q} + i_3 \\
&= x_i.
\end{align*}
Therefore, $0 \leq x_1 \leq \dots x_{n-1}$. 

\noindent \textbf{Part C.2: showing that for all $1\leq i,j < n$, we have $x_{T_{\cF}(i,j)} \leq x_i + x_j$.}
Fix $i,j$ and let $k = T_{\cF}(i,j)$.

\noindent \textbf{Case 1: $1 \leq i < q$ and $1 \leq j < q$.}
Then inspection shows that $k = \min(q-1, i+j)$. Thus
\[
	x_k = \frac{k}{2q} \leq \frac{i}{2q} + \frac{j}{2q} = x_i + x_j.
\]
\noindent \textbf{Case 2: $1 \leq i < q$ and $q \leq j < pq$.}
Then because $F_{pq-1} = \QQ(e_1,e_2)$ (see \Cref{eqn:p2q-fpq}), we have $j < k < pq$. Thus
\[
	x_k = 1 \leq \frac{i}{2q} + 1 = x_i + x_j.
\]
\noindent \textbf{Case 3: $1 \leq i < q$ and $pq \leq j < n$.}
Write $j = j_1 + j_2p + j_3pq$ in mixed radix notation with respect to $(p,q,p)$. Then $j < k \leq j_1 + \min(q-1,i+j_2)p + j_3pq$. Then
\begin{align*}
x_k &\leq x_{j_1 + \min(q-1,i+j_2)p + j_3pq} \\
&= j_1 \frac{1}{4pq} + \min(q-1,i+j_2) \frac{1}{2q} + j_3 \\
&\leq \frac{i}{2q} + j_1 \frac{1}{4pq} + j_2 \frac{1}{2q} + j_3 \\
&= x_i + x_j.
\end{align*}

\noindent \textbf{Case $4$: $q \leq i < pq$ and $q \leq j < pq$.}
Then as $F_{pq-1} = \QQ( e_1,e_2 )$, we have $j < k < pq$. Thus
\[
	x_k = 1 \leq \frac{i}{2q} + 1 = x_i + x_j.
\]
\noindent \textbf{Case $5$: $q \leq i < pq$ and $pq \leq j < n$.}
Write $j = j_1 + j_2p + j_3pq$ in mixed radix notation with respect to $(p,q,p)$. Because $F_{pq-1} = \QQ( e_1,e_2)$, then $j < k \leq (p-1) + (q-1)p + j_3pq$. Then
\begin{align*}
x_k &\leq x_{(p-1) + (q-1)p + j_3pq} \\
& = (p-1) \frac{1}{4pq} + (q-1) \frac{1}{2q} + j_3 \\
& \leq 1 + j_3 \\
&= 1 + j_1 \frac{1}{4pq} + j_2 \frac{1}{2q} + j_3  \\
&= x_i + x_j.
\end{align*}

\noindent \textbf{Case $6$: $pq \leq i < n$.}
Write
\begin{align*}
i &= i_1 + i_2p + i_3pq \\
j &= j_1 + j_2p + j_3pq
\end{align*}
in mixed radix notation with respect to $(p,q,p)$. Recall that for $i \geq pq$, we have $v_i = e_2^{i_1}e_1^{i_2}e_3^{i_3}$. Thus, $j < k \leq \min(p-1,i_1+j_1) + \min(q-1,i_2+j_2)p + (i_3 + j_3)pq$.  Then
\begin{align*}
x_k &\leq x_{\min(p-1,i_1+j_1) + \min(q-1,i_2+j_2)p + (i_3 + j_3)pq} \\
& = \min(p-1,i_1+j_1) \frac{1}{4pq} + \min(q-1,i_2+j_2)\frac{1}{2q} + (i_3 + j_3) \\
&\leq i_1 \frac{1}{4pq} + i_2 \frac{1}{2q} + i_3 + j_1 \frac{1}{4pq} + j_2 \frac{1}{2q} + j_3 \\
&=x_i + x_j.
\end{align*}

Thus, for all $1 \leq i,j < n$, we have if $x_{T_{\cF}(i,j)} \leq x_i + x_j$. 

\noindent \textbf{Part D: showing that $\mathbf{x} \notin \Cone(\Len_{(p,p,q)} \cup \Len_{(p,q,p)} \cup \Len_{(q,p,p)})$.} 
By definition 
\[
	\Len_{(p,p,q)} \cup \Len_{(p,q,p)} \subseteq \{\mathbf{x} \in \RR^{p^2q-1} : x_{q} \leq x_1 + x_{q-1} \}
\]
and
\[
	\Len_{(q,p,p)} \subseteq \{\mathbf{x} \in \RR^{p^2q-1} : x_{2pq+p} \leq x_{pq+1} + x_{pq + p - 1}\}.
\]
This implies that: 
\[
	\mathbf{x} \notin \Cone(\Len_{(p,p,q)} \cup \Len_{(p,q,p)} \cup \Len_{(q,p,p)}),
\]
which completes our proof.
\end{proof}

\begin{lemma}
\label{lem:helper-helper}
Let $q$ be an odd prime. For $a \in \ZZ/q\ZZ$, we have
\[
a\bigg\{\frac{q+1}{2},\dots,q-1\bigg\} = \bigg\{\frac{q+1}{2},\dots,q-1\bigg\} \Mod{q}
\]
if and only if $a \equiv 1 \Mod{q}$.
\end{lemma}
\begin{proof}
The statement 
\[
a\bigg\{\frac{q+1}{2},\dots,q-1\bigg\} = \bigg\{\frac{q+1}{2},\dots,q-1\bigg\} \Mod{q}
\]
is equivalent to the statement 
\[
a\bigg\{1,\dots,\frac{q-1}{2}\bigg \} = \bigg\{1,\dots,\frac{q-1}{2}\bigg \} \Mod{q}.
\]
If $q = 3$, it is clear that $a \equiv 1 \Mod{q}$; assume $q \neq 3$, and hence $q \geq 5$. Then as $a(q-1) \equiv -a \in \{\frac{q+1}{2},\dots,q-1\} \Mod{q}$, we must have $a \in \{1,\dots,\frac{q-1}{2}\} \Mod{q}$. If $a \not\equiv 1 \Mod{q}$, then there exists $b \in \{1,\dots,\frac{q-1}{2}\} \Mod{q}$ such that $ab \in \{\frac{q+1}{2},\dots,q-1\} \Mod{q}$, which is a contradiction.
\end{proof}

\begin{lemma}
\label{lem:pqr-helper}
Let $p$, $q$, and $r$ be odd prime numbers such that $p < q \leq r$. There exists an integer $m$ such that 
\[
	q \leq m \leq \lfloor qr/2 \rfloor,
\]
the addition $pm + pm$ overflows modulo $q$, and the addition $m + m$ does not overflow modulo $q$ or modulo $r$.
\end{lemma}
\begin{proof}
If $q = r$ then let
\[
	m = q + \bigg \lfloor \frac{q}{p} \bigg \rfloor.
\] 
Note that $m\perc q =  \lfloor q/p \rfloor$, as $0 \leq \lfloor q/p \rfloor < q$. As $0 \leq 2 \lfloor \frac{q}{p} \rfloor < 2q/p \leq q$, we have $(2m)\perc q = 2\lfloor q/p \rfloor$ and thus $m\perc q + m\perc q = (2m)\perc q$, so $m + m$ does not overflow modulo $q$.

On the other hand,  $pm = pq + p\lfloor q/p\rfloor$ and $0 \leq p\lfloor q/p\rfloor < q$, so $(pm)\perc q = p\lfloor q/p\rfloor$. Because $p < q$, we have $q\perc p \leq q/2$. Therefore,
\begin{align*}
(pm)\perc q + (pm)\perc q &= 2p\bigg \lfloor \frac{q}{p} \bigg \rfloor \\
&= 2p\bigg (\frac{q}{p}-\frac{q\perc p}{p}\bigg) \\
&\geq  2p\bigg (\frac{q}{p}-\frac{q}{2p}\bigg) \\
&= q,
\end{align*}
so the addition $pm + pm$ overflows modulo $q$.

Because $p^{-1} \neq 1 \Mod{q}$, by \cref{lem:helper-helper} the set
\[
	\bigg\{1,\dots,\frac{q-1}{2}\bigg\} \bigcap p^{-1}\bigg\{\frac{q+1}{2},\dots,q-1\bigg\} \Mod{q}
\]
is nonempty. Choose an element $\ell \in \ZZ/q\ZZ$ contained in the set above. Observe that 
\[
	q(r-1)/2 + (q-1)/2 = \bigg \lfloor \frac{qr}{2} \bigg \rfloor.
\]
and let
\[
	q \leq \ell_1,\dots, \ell_{\frac{r+1}{2}} 
\]
be the lifts of $\ell$ to $[q, \lfloor qr/2 \rfloor]$. Because $q \neq r$, the lifts $\ell_1,\dots, \ell_{\frac{r+1}{2}}$ all have distinct values modulo $r$ by the Chinese remainder theorem. Thus, there exists $\ell_k$ such that $\ell_k \in \{0,\dots, \frac{r-1}{2}\} \Mod{r}$. Set $m = \ell_k$.

To see that the addition $pm + pm$ overflows modulo $q$, notice that $m \in p^{-1}\{\frac{q+1}{2},\dots,q-1\} \Mod{q}$, so $(pm)\perc q \geq \frac{q+1}{2}$, and hence 
\[
	(pm)\perc q + (pm)\perc q \geq q+1.
\]

To see that the addition $m + m$ does not overflow modulo $q$ or $r$, observe that $m \in \{1,\dots,\frac{q-1}{2}\} \Mod{q}$, hence
\[
	m\perc q + m\perc q < q.
\] 
Similarly, since $m \in \{1,\dots,\frac{r-1}{2}\} \Mod{r}$, we have
\[
	m\perc r + m\perc r < r.
\] 
\end{proof}

\begin{proposition}
\label{lem:pqr}

Let $n = pqr$ for primes $p$, $q$, and $r$ with $p < q \leq r$. Then there exists a flag $\cF$ of a degree $n$ number field such that 
\[
	P_{T_{\cF}} \not \subseteq \Len_{(p,q,r)} \cup \Len_{(p,r,q)} \cup \Len_{(q,p,r)} \cup \Len_{(q,r,p)} \cup \Len_{(r,p,q)} \cup \Len_{(r,q,p)}.
\]
\end{proposition}
\begin{proof}

It suffices to show that there exists a flag $\cF$ and a point
\[
	\mathbf{x} \in \Cone(P_{T_{\cF}}) \setminus \Cone( \Len_{(p,q,r)} \cup \Len_{(p,r,q)} \cup \Len_{(q,p,r)} \cup \Len_{(q,r,p)} \cup \Len_{(r,p,q)} \cup \Len_{(r,q,p)}).
\] 

\noindent \textbf{Part A: defining the flag $\cF$.} Choose $e_1,e_2,e_3 \in \overline{\QQ}$ such that:
\begin{itemize}
    \item $e_1$ has degree $p$;
    \item $e_2$ has degree $q$;
    \item $e_3$ has degree $r$;
    \item and the compositum $K = \QQ(e_1,e_2,e_3)$ has degree $p^2q$.
\end{itemize}

Define a basis ${1=v_0,\dots,v_{n-1}}$ of $K$ as follows. For $0 \leq i < p$, set $v_i \coloneqq e_1^i$. For $p \leq i < n$ and $1 \leq i' < n$, write $i' = i'_1 + i'_2q + i_3'pq$ in mixed radix notation with respect to $(q,p,r)$. Inductively define $v_i$ as follows. Choose $i'$ minimal such that $e_2^{i'_1}e_1^{i'_2}e_3^{i'_3} \notin \{v_0,\dots,v_{i-1}\}$. Set $v_i \coloneqq e_2^{i'_1}e_1^{i'_2}e_3^{i'_3}$. Observe that for $i \geq pq$, we have $i = i'$. Define a flag $\cF = \{F_i\}_{i \in [n]}$ by $F_i \coloneqq \QQ\la v_0,\dots,v_i \ra$.

By \cref{lem:pqr-helper}, there exists an integer $m$ such that 
\[
	q \leq m \leq \lfloor qr/2 \rfloor,
\]
the addition $pm + pm$ overflows modulo $q$, and the addition $m + m$ does not overflow modulo $q$ or modulo $r$. Moreover, 
\[
	2pm \leq 2p \lfloor qr/2 \rfloor < pqr.
\]

Write $pm = (pm)_1 + (pm)_2 q + (pm)_3 pq$ in mixed radix notation with respect to $(q,p,r)$. 

\noindent \textbf{Part B: explicit description of $\cF$.}
We have:
\begin{align*}
F_1 &= \QQ \la 1, e_1 \ra \\
F_{p-1} &= \QQ( e_1 ) \\
F_{pm} &= \QQ\la \{e_2^{i_1}e_1^{i_2}e_3^{i_3} : i_1 + i_2q + i_3pq \leq pm, \; 0 \leq i_1 < q, \; 0 \leq i_2 < p, \; 0 \leq i_3\}\ra \\
F_{2pm-1} &= \QQ\la \{e_2^{i_1}e_1^{i_2}e_3^{i_3} : i_1 + i_2q + i_3pq \leq 2pm-1, \; 0 \leq i_1 < q, \; 0 \leq i_2 < p, \; 0 \leq i_3 \}\ra.
\end{align*}
We have 
\[
F_1F_{p-1}  = F_{p-1}.
\]
Moreover, because the addition $pm + pm$ overflows modulo $q$, we have $(pm)_1 + (pm)_1 \geq q$. We have that
\begin{align*}
F_{pm}F_{pm} \subseteq \QQ \la \{e_2^{i_1}e_1^{i_2}e_3^{i_3} : \; &i_1 + i_2q + i_3pq \leq (q-1) + \min(p-1,2(pm)_2) q + 2(pm)_3 pq, \\
&\; 0 \leq i_1 < q, \; 0 \leq i_2 < p, \; 0 \leq i_3\} \ra.
\end{align*}

Because $i_1 + i_2q + i_3pq \leq (q-1) + \min(p-1,2(pm)_2) q + 2(pm)_3 pq \leq 2pm-1$, we have
\[
	F_{pm}F_{pm} \subseteq F_{2pm-1}.
\]

\noindent \textbf{Part C: showing there exists $\mathbf{x} \in P_{T_{\cF}}$ such that $x_p > x_1 + x_{p-1}$ and $x_{2pm} > 2x_{pm}$.} For $1 \leq i < p$ set $x_i \coloneqq \frac{i}{2p}$. For $p \leq i < pq$ set $x_i \coloneqq 1$. For $pq \leq i < pqr$ write $i = i_1 + i_2q + i_3pq$ in mixed radix notation with respect to $(q,p,r)$ and set $x_i \coloneqq i_1 \frac{1}{4pq} + i_2 \frac{1}{2p} + i_3$. 

\noindent \textbf{Part C.1: showing that $0 \leq x_1 \leq \dots \leq x_{n-1}$.}
If $1 \leq i < p-1$, then
\[
	x_{i + 1} = \frac{i+1}{2p} \geq \frac{i}{2p} = x_i.
\]
Note also that
\[
	x_p = 1 \geq \frac{p-1}{2p} = x_{p-1}.
\]
If $p \leq i < pq-1$,
\[
	x_{i+1} = 1 = x_i.
\]
If $pq \leq i < n-1$ then write $i+1 = (i+1)_1 + (i+1)_2q + (i+1)pq$ in mixed radix notation with respect to $(q,p,r)$. If $i_1=q-1$ and $i_2 = p-1$ then $(i+1)_1 = (i+1)_2 = 0$ and $(i+1)_3 = i_3 + 1$. Then
\begin{align*}
x_{i+1} &= (i+1)_1 \frac{1}{4pq} + (i+1)_2 \frac{1}{2p} + (i+1)_3 \\
&= i_3 + 1 \\
&> i_1 \frac{1}{4pq} + i_2 \frac{1}{2p} + i_3 \\
&= x_{i}.
\end{align*}

Instead, if $i_1=q-1$ and $i_2 \neq p-1$, then $(i+1)_1 = 0$ and $(i+1)_2=i_2+1$ and $(i+1)_3 = i_3$. Then
\begin{align*}
x_{i+1} &= (i+1)_1 \frac{1}{4pq} + (i+1)_2 \frac{1}{2p} + (i+1)_3 \\
&= (i_2+1) \frac{1}{2p} + i_3 \\
&> i_1 \frac{1}{4pq} + i_2 \frac{1}{2p} + i_3 \\
&= x_{i}.
\end{align*}

Finally, if $i \neq q-1$ then $(i+1)_1 = i_1+1$ and $(i+1)_2=i_2$ and $(i+1)_3 = i_3$. Then
\begin{align*}
x_{i+1} &= (i+1)_1 \frac{1}{4pq} + (i+1)_2 \frac{1}{2p} + (i+1)_3 \\
&= (i_1+1) \frac{1}{4pq} + i_2 \frac{1}{2p} + i_3 \\
&> i_1 \frac{1}{4pq} + i_2 \frac{1}{2p} + i_3 \\
&= x_{i}.
\end{align*}

Therefore, $0 \leq x_1 \leq \dots \leq x_{n-1}$.

\noindent \textbf{Part C.2: showing that for all $1 \leq i,j < n$, we have $x_{T_{\cF}(i,j)} \leq x_i + x_j$.}

\noindent \textbf{Case $1$: $1 \leq i < p$ and $1 \leq j < p$.}
Then $k = \min(p-1, i+j)$. Then
\[
	x_k = \frac{k}{2p} \leq \frac{i}{2p} + \frac{j}{2p} = x_i + x_j.
\]
\noindent \textbf{Case $2$: $1 \leq i < p$ and $p \leq j < pq$.}
Then as $F_{pq-1} = \QQ( e_1,e_2)$, we have $j < k < pq$. Thus
\[
	x_k = 1 \leq \frac{i}{2p} + 1 = x_i + x_j.
\]
\noindent \textbf{Case $3$: $1 \leq i < p$ and $pq \leq j < n$.}
Write $j = j_1 + j_2q + j_3pq$ in mixed radix notation with respect to $(q,p,r)$. Then $j < k \leq j_1 + \min(p-1,i+j_2)q + j_3pq$. Then
\begin{align*}
x_k &\leq x_{j_1 + \min(p-1,i+j_2)q + j_3pq} \\
&= j_1 \frac{1}{4pq} + \min(p-1,i+j_2) \frac{1}{2p} + j_3 \\
&\leq \frac{i}{2p} + j_1 \frac{1}{4pq} + j_2 \frac{1}{2p} + j_3 \\
&= x_i + x_j.
\end{align*}

\noindent \textbf{Case $4$: $p \leq i < pq$ and $p \leq j < pq$.}
Then as $F_{pq-1} = \QQ( e_1,e_2)$, we have $j < k < pq$. Thus
\[
	x_k = 1 \leq \frac{i}{2q} + 1 = x_i + x_j.
\]
\noindent \textbf{Case $5$: $p \leq i < pq$ and $pq \leq j < n$.}
Write $j = j_1 + j_2q + j_3pq$ in mixed radix notation with respect to $(q,p,r)$. Because $F_{pq-1} = \QQ( e_1,e_2)$, we have $j < k \leq (q-1) + (p-1)q + j_3pq$. Then
\begin{align*}
x_k &\leq x_{(q-1) + (p-1)q + j_3pq} \\
& = (q-1) \frac{1}{4pq} + (p-1) \frac{1}{2p} + j_3 \\
& \leq 1 + j_3 \\
&= 1 + j_1 \frac{1}{4pq} + j_2 \frac{1}{2p} + j_3  \\
&= x_i + x_j.
\end{align*}

\noindent \textbf{Case $6$: $pq \leq i < n$.}
Write 
\begin{align*}
i &= i_1 + i_2q + i_3pq \\
j &= j_1 + j_2q + j_3pq
\end{align*}
in mixed radix notation with respect to $(q,p,r)$. Recall that for $i \geq pq$, we have $v_i = e_2^{i_1}e_1^{i_2}e_3^{i_3}$. Thus, $j < k \leq \min(q-1,i_1+j_1) + \min(p-1,i_2+j_2)q + (i_3 + j_3)pq$.  Then
\begin{align*}
x_k &\leq x_{\min(q-1,i_1+j_1) + \min(p-1,i_2+j_2)p + (i_3 + j_3)pq} \\
& = \min(q-1,i_1+j_1) \frac{1}{4pq} + \min(p-1,i_2+j_2)\frac{1}{2p} + (i_3 + j_3) \\
&\leq (i_1 \frac{1}{4pq} + i_2 \frac{1}{2p} + i_3) + (j_1 \frac{1}{4pq} + j_2 \frac{1}{2p} + j_3) \\
&=x_i + x_j.
\end{align*}

Thus, for any integers $1 \leq i,j < n$, we have $x_{T_{\cF}(i,j)} \leq x_i + x_j$. 

\noindent \textbf{Part C.3: showing that $x_p > x_1 + x_{p-1}$ and $x_{2pm} > 2x_{pm}$.}
Moreover, we have that
\[
	x_p = 1 > \frac{1}{2p} + \frac{p-1}{2p} = x_1 + x_{p-1}.
\]
Write: 
\begin{align*}
pm &= (pm)_1 + (pm)_2q + (pm)_3pq \\
2pm &= (2pm)_1 + (2pm)_2q + (2pm)_3pq
\end{align*}
in mixed radix notation with respect to $(q,p,r)$ and recall that $pm + pm$ overflows modulo $q$. Therefore, either $(2pm)_3 = 2(pm)_3$ and $(2pm)_2 = 2(pm)_2 + 1$, or $(2pm)_3 = 2(pm)_3 + 1$. If $(2pm)_3 = 2(pm)_3$ and $(2pm)_2 > 2(pm)_2+1$ then
\begin{align*}
x_{2pm} &= (2pm)_1 \frac{1}{4pq} + (2pm)_2\frac{1}{2p} + (2pm)_3 \\
&\geq (2(pm)_2+1)\frac{1}{2p} + 2(pm)_3 \\
&> 2(pm)_1\frac{1}{4pq} + 2(pm)_2\frac{1}{2p} + 2(pm)_3 \\
&= 2x_{pm}.
\end{align*}
Otherwise, if $(2pm)_3 = 2(pm)_3 + 1$, then 
\begin{align*}
x_{2pm} &= (2pm)_1 \frac{1}{4pq} + (2pm)_2\frac{1}{2p} + (2pm)_3 \\
&\geq 2(pm)_3 + 1 \\
&> 2(pm)_1\frac{1}{4pq} + 2(pm)_2\frac{1}{2p} + 2(pm)_3 \\
&= 2x_{pm}.
\end{align*}

\noindent \textbf{Part D: showing that $\mathbf{x} \notin \Cone(\Len_{(q,p,r)} \cup \Len_{(q,r,p)} \cup \Len_{(r,q,p)} \cup \Len_{(r,p,q)} \cup \Len_{(p,q,r)} \cup \Len_{(p,r,q)})$.}
Note that 
\[
	\Len_{(q,p,r)} \cup \Len_{(q,r,p)} \cup \Len_{(r,q,p)} \cup \Len_{(r,p,q)} \subseteq \{\mathbf{x} \in \RR^{n-1} : x_{p} \leq x_1 + x_{p-1} \}
\]
and
\[
	\Len_{(p,q,r)} \cup \Len_{(p,r,q)} \subseteq \{\mathbf{x} \in \RR^{n-1} : x_{2pm} \leq 2x_{pm}\}.
\]
We have $x_p > x_1 + x_{p-1}$ and $x_{2pm} > 2x_{pm}$, and thus our proof is complete. 
\end{proof}

\begin{proposition}
\label{lem:4p}
Let $n = 4p$ for $p$ a prime not equal to $2$ or $3$. Then there exists a flag $\cF$ of a degree $n$ number field such that 
\[
	P_{T_{\cF}} \not \subseteq \Len_{(2,2,p)} \cup \Len_{(2,p,2)} \cup \Len_{(p,2,2)}.
\]
\end{proposition}
\begin{proof}
It suffices to show that there exists a flag $\cF$ and a point 
\[
	\mathbf{x} \in \Cone(P_{T_{\cF}}) \setminus \Len_{(2,2,p)} \cup \Len_{(2,p,2)} \cup \Len_{(p,2,2)}.
\]

\noindent \textbf{Part A: defining the flag $\cF$.}
Choose $e_1,e_2,e_3 \in \overline{\QQ}$ such that:
\begin{itemize}
\item $e_1$ has degree $p$;
\item the element $e_2$ has degree $2$;
\item  the element $e_3$ has degree $2$;
\item and the compositum $\QQ(e_1,e_2,e_3)$ has degree $4p$. 
\end{itemize}
 
Define a basis $\{ v_0,\dots,v_{4p-1}\}$ of $K$ as follows. Set 
\begin{align*}
v_0 &\coloneqq 1 \\
v_1 &\coloneqq e_1 \\
v_2 &\coloneqq e_2 \\
v_3 &\coloneqq e_2 e_1 \\
v_4 &\coloneqq e_1^2 \\
v_5 &\coloneqq e_2 e_1^2 \\
v_6 &\coloneqq e_1^3 \\
v_7 &\coloneqq e_2e_1^3.
\end{align*}
For $8 \leq i < n$ and $1 \leq i' < n$, write $i' = i'_1 + i'_2p + i_3'2p$ in mixed radix notation with respect to $(p,2,2)$. Inductively define $v_i$ as follows. Choose $i'$ minimal such that $e_1^{i'_1}e_2^{i'_2}e_3^{i'_3} \notin \{v_0,\dots,v_{i-1}\}$. Set $v_i = e_1^{i'_1}e_2^{i'_2}e_3^{i'_3}$. Observe that for $i \geq 2p$, we have $i = i'$. Let $\cF$ be the corresponding flag.

\noindent \textbf{Part B: explicit description of $\cF$.}
Note that:
\begin{align*}
F_1 &= \QQ\la  1, e_1 \ra \\
F_3 &= \QQ( e_2) \QQ \la 1,e_1 \ra\\
F_5 &= \QQ( e_2)\QQ\la  1, e_1, e_1^2 \ra \\
F_7 &= \QQ( e_2)\QQ \la 1, e_1, e_1^2, e_1^3 \ra \\
F_{3p-1} &= \QQ( e_1 )  \QQ \la 1, e_2, e_3 \ra.
\end{align*}
Therefore
\[
	F_1F_{3p-1}  = F_{3p-1}
\]
\[
F_3F_3 = F_5
\]
\[
F_3F_5 = F_7.
\]

\noindent \textbf{Part C: defining $\mathbf{x}$ when $p = 5$.}
Suppose $p=5$. Then let $\mathbf{x} \in \RR^{19}$ be as follows:
\begin{align*}
x_1 &\coloneqq 1 \\
x_2,x &\coloneqq 1.4 \\
x_4,x_5 &\coloneqq 2 \\
x_6,x_7 &\coloneqq 3 \\
x_8,\dots,x_{14} &\coloneqq 4 \\
x_{15},\dots,x_{19} &\coloneqq 5.1.
\end{align*}

\noindent \textbf{Part D: showing that $\mathbf{x} \in \Cone(P_{T_{\cF}})$ and $x_8 > x_5 + x_3$ and $x_{15} > x_1 + x_{14}$ when $p = 5$.}
It is clear that $0 \leq x_1 \leq \dots \leq x_{n-1}$ and $x_8 > x_5 + x_3$ and $x_{15} > x_1 + x_{14}$. We now show that for all $1\leq i,j < n$, we have $x_{T_{\cF}(i,j)} \leq x_i + x_j$. Fix $i,j$ and let $k = T_{\cF}(i,j)$.

\noindent \textbf{Case 1: $i = 1$.} If $j = 1$ then $k = 4$ and 
\[
	x_4 = 2 \leq 2x_1.
\]
Observe that if $j = 2,3$ then $k = 4,5$, and
\[
	x_k = 2 \leq 1.4 + 1 \leq x_1 + x_j.
\]
If $j = 4,5$, then $k = 6,7$, and 
\[
	x_k = 3 \leq 2 + 1 = x_1 + x_j.
\]
If $j = 6,7$, then $k \leq 10$, and
\[
	x_k \leq 4 \leq 3 + 1 = x_1 + x_j.
\]
If $8 \leq j < 15$, then as $v_1 = e_1$ and $F_{14} = \la e_1 \ra \{ 1, e_2, e_3 \}$, we have that $j < k < 14$. Thus
\[
	x_k \leq 4 \leq 4 + 1 = x_1 + x_j.
\] 
If $15\leq j < n$ then
\[
	x_k \leq 5.1 \leq 5.1 + 1 \leq x_1 + x_j.
\]

\noindent \textbf{Case 2: $i = 2,3$.}
If $j = 2,3$, then $k = 4,5$ so 
\[
	x_k = 2 \leq 1.4 + 1.4 = x_i + x_j.
\]
If $j = 4,5$ then $k = 6,7$ so
\[
	x_k = 3 \leq 2 + 1.4 = x_i + x_j.
\]
If $j = 6,7$ then $k < 15$ so
\[
	x_k \leq 4 \leq 3 + 1.4 \leq x_i + x_j.
\]
If $8 \leq j < n$ then
\[
	x_k \leq 5.1 \leq 4 + 1.4 \leq x_i + x_j.
\]

\noindent \textbf{Case 3: $i = 4,5$.}
If $j < 10$ then because $F_9 = \QQ( e_1,e_2 )$ we have $k < 10$. Thus
\[
	x_k \leq 4 \leq 2 + 2 \leq x_i + x_j. 
\]
If $10 \leq j < n$ then
\[
	x_k \leq 5.1 \leq 4 + 2 \leq x_i + x_j.
\]

\noindent \textbf{Case 4: $i \geq 6$.}
Then 
\[
	x_k \leq 5.1 \leq 3 + 3 \leq x_i + x_j.
\]

\noindent \textbf{Part E: finishing the proof when $p = 5$.}
Because 
\[
	\Len_{(2,2,5)} \cup \Len_{(2,5,2)} \subset \{\mathbf{x} \in \RR^{19} : x_{15} \leq x_1 + x_{14} \}
\]
and
\[
	\Len_{(5,2,2)} \subset \{\mathbf{x} \in \RR^{19} : x_8 \leq x_5 + x_3 \},
\]
we have that 
\[
	\mathbf{x} \notin \Cone(\Len_{(2,2,5)} \cup \Len_{(2,5,2)} \cup \Len_{(5,2,2)})
\]
Hence, the proof is complete for $p = 5$.

\noindent \textbf{Part F: defining $\mathbf{x}$ when $p \neq 5$.}
If $p \neq 5$, let $\mathbf{x} \in \RR^{4p-1}$ be as follows.
\begin{align*}
x_1 &\coloneqq 1 \\
x_2, x_3 &\coloneqq 1.4 \\
x_4, x_5 &\coloneqq 2 \\
x_6, \dots, x_{3p-1} &\coloneqq 2.9 \\
x_{3p}, \dots, x_{4p-1} &\coloneqq 4.
\end{align*}

\noindent \textbf{Part G: showing that $\mathbf{x} \in \Cone(P_{T_{\cF}})$ and $x_6 > x_3 + x_3$ and $x_{3p} > x_1 + x_{3p-1}$.}
It is clear that $0\leq x_1 \leq \dots \leq x_{n-1}$ and $x_6 > x_3 + x_3$ and $x_{3p} > x_1 + x_{3p-1}$. We now show that for all integers $1 \leq i,j < n$, we have $x_{T_{\cF}(i,j)} \leq x_i + x_j$. Fix $i,j$ and let $k = T_{\cF}(i,j)$.

\noindent \textbf{Case 1: $i = 1$.}
If $j = 1$, then $k = 4$ and
\[
	x_4 = 2 \leq 2x_1.
\]
Observe that if $j = 2,3$, then $k = 4,5$, and
\[
	x_k \leq x_5 = 2 \leq 1.4 + 1 \leq x_1 + x_j.
\]
If $j = 4,5$, then $k = 6,7$, and
\[
	x_k \leq x_5 = 2.9 \leq 2 + 1 \leq x_1 + x_j.
\]
If $6 \leq j < 3p$, then as $v_1 = e_1$ and $F_{3p-1} = \QQ( e_1 ) \QQ \la 1, e_2, e_3 \ra$, we have $j < k < 3p$. Thus
\[
	x_k = 2.9 \leq 2.9 + 1 = x_j + x_1.
\]
If $j \geq 3p$, then
\[
	x_k = 4 \leq 4 + 1 \leq x_1 + x_j.
\]
\textbf{Case 2: $i = 2,3$.} 
If $j = 2,3$, then $k \leq 5$ and
\[
	x_k = 2 \leq 1.4 + 1.4 = x_i + x_j.
\]
If $j = 4,5$ then $k \leq 7$ and
\[
	x_k \leq 2.9 \leq 2 + 1.4 = x_i + x_j.
\]
If $j \geq 6$, then
\[
	x_k \leq 4 \leq 2.9 + 1.4 = x_i + x_j.
\]
\textbf{Case 3: $i \geq 4$.}
Then we have
\[
	x_k \leq 4 \leq 2 + 2 \leq x_i + x_j.
\]

\noindent \textbf{Part H: finishing the  proof when $p \neq 5$.}
By definition, we have
\[
	\Len_{(2,2,p)} \cup \Len_{(2,p,2)} \subset \{\mathbf{x} \in \RR^{4p-1} : x_{3p} \leq x_1 + x_{3p-1} \}
\]
and because $p \geq 7$, we have
\[
	\Len_{(p,2,2)} \subset \{\mathbf{x} \in \RR^{4p-1} : x_6 \leq x_3 + x_3 \},
\]
Therefore, $\mathbf{x} \notin \Cone(\Len_{(2,2,p)} \cup \Len_{(2,p,2)} \cup \Len_{(p,2,2)})$.
\end{proof}

\begin{proposition}
\label{lem:24}
Let $n = 24$. Then there exists a flag $\cF$ of a degree $n$ number field such that 
\[
	P_{T_{\cF}} \not \subseteq \Len_{(2,2,2,3)} \cup \Len_{(2,2,3,2)} \cup \Len_{(2,3,2,2)} \cup \Len_{(3,2,2,2)}
\]
\end{proposition}
\begin{proof}
As usual, it suffices to show that there exists a flag $\cF$ and a point
\[
	\mathbf{x} \in P_{T_{\cF}} \setminus  \Len_{(2,2,2,3)} \cup \Len_{(2,2,3,2)} \cup \Len_{(2,3,2,2)} \cup \Len_{(3,2,2,2)}.
\]

\noindent \textbf{Part A: defining the flag $\cF$.}
Choose elements $e_1,e_2,e_3,e_4 \in \overline{\QQ}$ such that:
\begin{itemize}
\item $e_1,e_2,e_3$ have degree $2$;
\item $e_4$ has degree $3$;
\item and the compositum $K= \QQ(e_1,e_2,e_3,e_4)$ has degree $24$.
\end{itemize}
Define a basis $\{v_0,\dots,v_{23}\}$ of $K$ via the formulae:
\begin{align*}
v_0 & \coloneqq 1 \\
v_1 & \coloneqq e_4 \\
v_2 & \coloneqq e^2_4 \\
v_3 & \coloneqq e_1 \\
v_4 & \coloneqq e_1e_4 \\
v_5 & \coloneqq e_1e_4^2 \\
v_6 & \coloneqq e_2 \\
v_7 & \coloneqq e_1e_2 \\
v_8 & \coloneqq e_4e_2 \\
v_9 & \coloneqq e_1e_4e_2 \\
v_{10} & \coloneqq e^2_4e_2 \\
v_{11} & \coloneqq e_1e^2_4e_2 \\
v_{12} & \coloneqq e_3 \\
v_{13} & \coloneqq e_1e_3 \\
v_{14} & \coloneqq e_4e_3 \\
v_{15} & \coloneqq e_1e_4e_3 \\
v_{16} & \coloneqq e^2_4e_3 \\
v_{17} & \coloneqq e_1e^2_4e_3 \\
v_{18} & \coloneqq e_2e_3 \\
v_{19} & \coloneqq e_1e_2e_3 \\
v_{20} & \coloneqq e_4e_2e_3 \\
v_{21} & \coloneqq e_1e_4e_2e_3 \\
v_{22} & \coloneqq e^2_4e_2e_3 \\
v_{23} & \coloneqq e_1e^2_4e_2e_3. \\
\end{align*}

\noindent Let $\cF$ be the associated flag.  

\noindent \textbf{Part C: showing that there exists $\mathbf{x} \in P_{T_{\cF}}$ with $x_3 > x_2 + x_1$ and $x_{20} > x_7 + x_{13}$.}
This can be checked explicitly using the computer algebra system Magma.

\noindent \textbf{Part D: finishing the proof.}
By definition, we have that: 
\[
	\Len_{(2,2,2,3)} \cup \Len_{(2,2,3,2)} \cup \Len_{(2,3,2,2)} \subset \{\mathbf{x} \in \RR^{23} : x_{3} \leq x_2 + x_1 \}
\]
and 
\[
	\Len_{(3,2,2,2)} \subset \{\mathbf{x} \in \RR^{23} : x_{20} \leq x_7 + x_{13} \}.
\]
Therefore, 
\[
	\mathbf{x} \notin \Len_{(2,2,2,3)} \cup \Len_{(2,2,3,2)} \cup \Len_{(2,3,2,2)}  \cup \Len_{(3,2,2,2)}.
\]
\end{proof}

\section{Bounds on scrollar invariants of curves}
\label{sec:fn-field}

In this section we switch focus and prove bounds on scrollar invariants of curves. Namely, we prove \Cref{thm:bound-subfield-fn-field} and \Cref{thm:containment-fn-field}. Let $k$ be a field and let $\pi \colon C \rightarrow \PP_1^k$ be a finite morphism from a smooth projective geometrically irreducible curve over $k$ to $\PP^1_k$. For conciseness, let $e_j^i \coloneqq e_j(\cL_i)$ denote the $j$th scrollar invariant of $\cL_i$. 

\begin{observation}
Consider $\cO_{\PP^1}$-algebra structure on the locally free module $\pi_* \cL_i = \cO_{\PP^1}(-e^i_0) \oplus \cO_{\PP^1}(-e^i_1) \oplus \cdots \oplus \cO_{\PP^1}(-e^i_{n-1})$. The map $\cL_1 \otimes \cL_2 \rightarrow \cL_3$ induces a map
\[
    \pi_* \cL_1 \otimes \pi_* \cL_2 \rightarrow \pi_* \cO_3. 
\]
Under the product structure in this sheaf of algebras, the product of the $i$th summand and the $j$th summand, decomposed again into summands, must be zero in any summand $\cO(-e^3_k)$ where $e^3_k > e^1_i + e^2_j$, as $\Hom(\cO(-e^1_i) \otimes \cO(-e^2_j), \cO(-e^3_k)) = 0$ in that case.  
\end{observation}

Choose a point $\infty \in \PP^1$, and choose a coordinate $t$ on $\AAA^1 = \PP^1 \setminus \{ \infty \}$, i.e., an isomorphism $\PP^1 \setminus \{ \infty \}  \cong \Spec k[t]$.  Then the splitting $\pi_* \cL_i = \cO_{\PP^1}(-e^i_0) \oplus \cdots \oplus \cO_{\PP^1}(-e^i_{n-1})$ induces a splitting of the $k[t]$-algebra
$\Gamma(\AAA^1, \pi_*(\cL_i))$ into
$\Gamma(\AAA^1, \pi_*(\cO_C)) = k[t] \oplus k[t] x_1 \oplus \cdots \oplus k[t] x_{n-1}$
(as a $k[t]$-module); here we have chosen a generator  $x_j$ of the $j$th summand  of $\Gamma(\AAA^1, \pi_*(\cL_i))$.
Considered as a rational section of $\pi_* \cL_i$, $x_j$ has a pole at $\infty$ of order $e^i_j$.  Now $1=x_0, x_1$, \dots, $x_{n-1}$ form a basis for $K(C)$ as a $K(\PP^1)$-vector space (where $K(\cdot)$ indicates the function field). 

\begin{proof}[Proof of \Cref{thm:bound-subfield-fn-field}]
Let $x_0,\dots,x_{n-1}$ (resp. $y_0,\dots,y_{n-1}$, $z_0,\dots,z_{n-1}$) be generators for $\pi_* \cL_1$ (resp. $\pi_* \cL_2$, $\pi_* \cL_3$). Let $I \coloneqq K(\PP^1)\langle x_0,\dots,x_i\rangle$ and $J \coloneqq K(\PP^1)\langle x_0,\dots,x_i\rangle$. If $\dim_{K(\PP^1)}IJ \geq i + j + 1$, then \Cref{prop:mult-bound} implies that
\[
    e_{i+j}(\cL_3) \leq e_i(\cL_1)+e_j(\cL_2),
\]
which is the desired conclusion. 

Now assume for the sake of contradiction that $\dim_{K(\PP^1)}IJ \leq i + j$ and set $m = \dim_{K(\PP^1)}\Stab(IJ)$.  The conclusion of \Cref{lem:overflow} states that $i_1 + j_1 \geq m$ and $m > 1$. Therefore,   
\[
    (i\%m) + (j\%m) \neq (i + j)\%m.
\]
However, this contradicts the assumptions of \Cref{thm:bound-subfield-fn-field} because $\Stab(IJ)$ is a field.
\end{proof}

\begin{proof}[Proof of \Cref{thm:containment-fn-field}]
Let $x_0,\dots,x_{n-1}$ be generators as above for $\pi_* \cO_C$. Let $\cF$ be the flag of $K(C)/K(\PP^1)$ given by $F_i = K(\PP^1)\la x_0,\dots,x_i\ra$ and let $T_{\cF}$ be the corresponding flag type. By \Cref{thm:flag-types}, there exists a tower type $\mfT$ such that $T_{\mfT} \leq T_{\cF}$. Let $i,j$ be integers such that $i+j$ does not overflow modulo $\mfT$. Then:
\begin{align*}
e_{i+j}(\cO_C) &=  e_{T_{\mfT}(i,j)}(\cO_C) && \text{ because $i+j = T_{\mfT}(i,j)$ by  \Cref{cor:explicit-flag-type}} \\
&\leq e_{T_{\cF}(i,j)}(\cO_C) && \text{because $T_{\mfT} \leq T$, so $T_{\mfT}(i,j) \leq T_{\cF}(i,j)$}\\
&\leq e_i(\cO_C)+ e_j(\cO_C) && \text{by \Cref{prop:flag-bound}}
\end{align*}
\end{proof}

\end{document}